\documentclass[oneside,english]{amsart}
\usepackage{pslatex}
\usepackage[T1]{fontenc}
\usepackage[latin1]{inputenc}
\usepackage{graphicx}
\usepackage{amssymb}

\makeatletter
 \theoremstyle{plain}    
 \newtheorem{thm}{Theorem}[section]
 \numberwithin{equation}{section} 
 \numberwithin{figure}{section} 
 \theoremstyle{plain}
 \theoremstyle{plain}    
 \newtheorem{lem}[thm]{Lemma} 

\usepackage{babel}
\makeatother
\begin{document}

\title{An efficient algorithm for accelerating the convergence of oscillatory
series, useful for computing the polylogarithm and Hurwitz zeta functions}

\author{Linas Vep\v{s}tas}

\date{12 October 2006, revised 18 March 2007}

\begin{abstract}
This paper sketches a technique for improving the rate of convergence
of a general oscillatory sequence, and then applies this series acceleration
algorithm to the polylogarithm and the Hurwitz zeta function. As such,
it may be taken as  an extension of the techniques given by Borwein's
{}``An efficient algorithm for computing the Riemann zeta function''\cite{Bor95,Bor00-A},
to more general series. The algorithm provides a rapid means of evaluating
$\mbox{Li}_{s}(z)$ for general values of complex $s$ and a kidney-shaped
region of complex $z$ values given by $\left|z^{2}/(z-1)\right|<4$.
By using the duplication formula and the inversion formula, the range
of convergence for the polylogarithm may be extended to the entire
complex $z$-plane, and so the algorithms described here allow for
the evaluation of the polylogarithm for all complex $s$ and $z$
values. 

Alternatively, the Hurwitz zeta can be very rapidly evaluated by means
of an Euler-Maclaurin series. The polylogarithm and the Hurwitz zeta
are related, in that two evaluations of the one can be used to obtain
a value of the other; thus, either algorithm can be used to evaluate
either function. The Euler-Maclaurin series is a clear performance
winner for the Hurwitz zeta, while the Borwein algorithm is superior
for evaluating the polylogarithm in the kidney-shaped region. Both
algorithms are superior to the simple Taylor's series or direct summation.

The primary, concrete result of this paper is an algorithm allows
the exploration of the Hurwitz zeta in the critical strip, where fast
algorithms are otherwise unavailable. 

A discussion of the monodromy group of the polylogarithm is included.
\end{abstract}



\maketitle

\section{Introduction}

This note sketches a technique for improving the rate of convergence
of a general oscillatory series, and then applies this technique to
the computation of the polylogarithm, and, in particular, to the Hurwitz
zeta function. It essentially generalizes an algorithm given by Peter
Borwein for computing the Riemann zeta function. The principle result
is an algorithm that efficiently obtains values of the Hurwitz zeta
in the critical strip $s=\sigma+i\tau$ with $0\le\sigma\le1$. 

The Hurwitz zeta function may be defined as 

\begin{equation}
\zeta(s,q)=\sum_{n=0}^{\infty}\frac{1}{(n+q)^{s}}\label{eq:Hurwitz zeta}\end{equation}
 This series may be directly summed for $\sigma>1$, although, for
high-precision work, convergence is annoyingly slow unless one has
at least $\sigma\gtrsim2$. A globally convergent series was given
by Helmut Hasse in 1930\cite{Has30}: \begin{equation}
\zeta(s,q)=\frac{1}{s-1}\sum_{n=0}^{\infty}\frac{1}{n+1}\sum_{k=0}^{n}(-1)^{k}\left(\begin{array}{c}
n\\
k\end{array}\right)(q+k)^{1-s}\label{eq:Hasse series}\end{equation}
 Although this series is technically convergent everywhere, in practice,
the convergence is abysmally slow on the critical strip. A not unreasonable
algorithm may be given by considering a Taylor's series in $q$. Expanding
about $q=0$, one obtains \begin{equation}
\zeta(s,q)=\frac{1}{q^{s}}+\sum_{n=0}^{\infty}\left(-q\right)^{n}\left(\begin{array}{c}
s+n-1\\
n\end{array}\right)\zeta(s+n)\label{eq:Hurwitz taylor}\end{equation}
 It is not hard to show that the sum on the right is convergent for
$\left|q\right|<1$. Alternately, an expansion may be made about $q=1/2$:
\[
\zeta\left(s,q+\frac{1}{2}\right)=\sum_{n=0}^{\infty}\left(-q\right)^{n}\left(\begin{array}{c}
s+n-1\\
n\end{array}\right)\left(2^{s+n}-1\right)\zeta(s+n)\]
 It may be shown that the above has a region of convergence $\left|q\right|<1/2$.
Either of these expansions provide a computable expression for the
Hurwitz zeta function that is convergent on the entire complex $s$-plane
(minus, of course, $s=1$, and taking the appropriate limit, via l'Hopital's
rule, when $s$ is a non-positive integer). The principal drawback
to the use of these sums for high-precision calculations is the need
for many high-precision evaluations of the Riemann zeta; the compute
time can grow exponentially as the number of required digits is increased,
especially when $q$ approaches the radius of convergence of the sums. 

It is the poor performance of these sums that motivates the development
of this paper. Since the Generalized Riemann Hypothesis can be phrased
in terms of the values of the Hurwitz zeta function on the critical
strip, it is of some interest to have a fast algorithm for computing
high-precision values of this function in this region. There seems
to be a paucity of work in this area. There is a discussion of fast
algorithms for Dirichlet L-functions in \cite{Sle04} (the author
has been unable to obtain a copy of this manuscript). 

The Hurwitz zeta may be expressed in terms of the polylogarithm\cite{Lew81}
(sometimes called the {}``\emph{fractional polylogarithm}'' to emphasize
that $s$ is not an integer): \begin{equation}
\mbox{Li}_{s}(z)=\sum_{n=1}^{\infty}\frac{z^{n}}{n^{s}}\label{eq:polylog definition}\end{equation}
 by means of Jonquière's identity\cite[Section 7.12.2]{Lew81}\cite{Jon1889}
\begin{equation}
\mbox{Li}_{s}\left(e^{2\pi iq}\right)+\left(-1\right)^{s}\mbox{Li}_{s}\left(e^{-2\pi iq}\right)=\frac{\left(2\pi i\right)^{s}}{\Gamma(s)}\zeta(1-s,q)\label{eq:poly-hur relation}\end{equation}
 and so one might search for good algorithms for polylogarithms. Wood\cite{Wood92}
provides a extensive review of the means of computing the polylogarithm,
but limits himself to real values of $s$. Crandall\cite{Cra06} discusses
evaluating the polylogarithm on the entire complex $z$-plane, but
limits himself to integer values of $s$. Thus, there appears to be
a similar paucity of general algorithms for the polylogarithm as well.
The series defining the polylogarithm may be directly evaluated when
$\left|z\right|<1$, although direct evaluation becomes quite slow
when $\left|z\right|\gtrsim0.95$ and $\sigma\lesssim2$. There do
not seem to be any published algorithms that may be applied on the
critical strip.

The primary effort of this paper is to take the algorithm given by
Borwein\cite{Bor95,Bor00-A}, which is a simple Padé-approximant type
algorithm, and generalize it to the polylogarithm. The result is a
relatively small finite sum that approximates the polylogarithm, and
whose error, or difference from the exact value, can be precisely
characterized and bounded. Increased precision is easily obtainable
by evaluating a slightly larger sum; one may obtain roughly $2N$
to $4N$ bits of precision by evaluating a sum of $N$ terms. The
sum may be readily evaluated for just about any value $s$ on the
complex $s$-plane. However, it is convergent only in a limited region
of the $z$-plane, and specifically, is only convergent when \[
\left|\frac{z^{2}}{z-1}\right|<4\]
 This is sufficient for evaluating the Hurwitz zeta for general complex
$s$ and real $0<q<1$. Unfortunately, there does not appear to be
any simple and efficient way of extending this result to the general
complex $z$-plane, at least not for general values of $s$. By using
duplication and reflection formulas, one may extend the algorithm
to the entire complex $z$-plane; however, this requires the independent
evaluation of Hurwitz zeta.

Although the Hurwitz zeta function may be computed by evaluating the
Taylor's series \ref{eq:Hurwitz taylor} directly, a superior approach
is to perform an Euler-Maclaurin summation (thanks to Oleksandr Pavlyk
for pointing this out\cite{Pav07}). The summation uses the standard,
textbook Euler-Maclaurin formula\cite[25,4,7]{A&S}, and is applied
to the function $f(x)=(x+q)^{-s}$. However, the application is {}``backwards''
from the usual sense: the integral of $f(x)$ is known (it is easily
evaluated analytically), and it is the series, which gives the Hurwitz
zeta, which is unknown. All of the other parts of Euler-Maclaurin
formula are easily evaluated. The result is an algorithm that is particularly
rapid for the Hurwitz zeta function. It outperforms direct evaluation
of the Taylor's series by orders of magnitude. It is faster then evaluating
\ref{eq:poly-hur relation} (which can be computed for real values
of $q$). 

The development of this paper proceeds as follows. First, a certain
specific integral representation is given for the polylogarithm. This
representation is such that a certain trick, here called {}``the
Borwein trick'', may be applied, to express the integral as a polynomial
plus a small error term. This is followed by a sketch of a generalization
of the trick to the convergence acceleration of general oscillatory
series. The next section selects a specific polynomial, and the error
term of the resulting approximation is precisely bounded. This is
followed by a very short review of the application of the Euler Maclaurin
summation. This is followed by a brief review of the duplication formula
for the Hurwitz zeta, and a short discussion of ways in which this
numerical algorithm may be tested for correctness. Measurements of
the performance of actual implementations of the various numerical
algorithms is provided. This is followed by a detailed derivation
of the monodromy group, and a discussion of Apostol's {}``periodic
zeta function''.

The algorithm has been implemented by the author using the Gnu Multiple
Precision arithmetic library\cite{GMP}, and is available on request,
under the terms of the LGPL license. The paper concludes with some
intriguing images of the polylogarithm and the Hurwitz zeta function
on the critical strip.

\section{The Polylogarithm}

The polylogarithm has a non-trivial analytic structure. For fixed
$s$, the principal sheet has a branch cut on the positive real axis,
for $1<z<+\infty$. The location of the branch-point at $z=1$, as
is always the case with branch points, is the cause of limited convergence
of analytic series in the complex $z$-plane. This is noted here,
as this branch point has a direct impact on the algorithm, preventing
convergence in its vicinity. Besides this, the polylogarithm has many
other interesting properties, which are not reviewed here.

The development of the algorithm requires the following integral representation.

\begin{lem}
The polylogarithm has the integral representation \[
\mbox{Li}_{s}\left(z\right)=\frac{z}{\Gamma(s)}\int_{0}^{1}\frac{\left|\log y\right|^{s-1}}{1-yz}\, dy\]

\end{lem}
\begin{proof}
This identity is easily obtained by inserting the integral representation
of the Gamma function:\begin{eqnarray*}
\mbox{Li}_{s}(z) & = & \frac{1}{\Gamma(s)}\,\sum_{n=1}^{\infty}\frac{z^{n}}{n^{s}}\int_{0}^{\infty}e^{-u}u^{s-1}du\\
 & = & \frac{1}{\Gamma(s)}\,\sum_{n=1}^{\infty}z^{n}\int_{0}^{\infty}e^{-nt}t^{s-1}dt\\
 & = & \frac{1}{\Gamma(s)}\,\int_{0}^{\infty}t^{s-1}\sum_{n=1}^{\infty}\left(ze^{-t}\right)^{n}\, dt\\
 & = & \frac{1}{\Gamma(s)}\,\int_{0}^{\infty}t^{s-1}\left[\frac{ze^{-t}}{1-ze^{-t}}\right]\, dt\end{eqnarray*}
 and then finally substituting $t=-\log y$ in the last integral.
Although this derivation casually interchanges the order of integration
and summation, one may appeal to general arguments about analytic
continuation to argue that the final result is generally valid, provided
one is careful to navigate about the branch point at $z=1$.
\end{proof}
A more sophisticated presentation of this theorem is given by Costin\cite[Thm. 1]{Cos07}.

\section{The Borwein trick}

The Borwein trick uses the integral form to find a simple, finite
sum that approximates the desired value arbitrarily well. The trick
consists of two steps. The first is to write the above integral in
the form \begin{eqnarray}
\xi(s,z) & = & \frac{1}{f(1/z)}\,\frac{z}{\Gamma(s)}\,\int_{0}^{1}\,\frac{f(y)\left|\log y\right|^{s-1}}{1-yz}\, dy\nonumber \\
 & = & \mbox{Li}_{s}\left(z\right)+\frac{1}{f(1/z)}\,\frac{z}{\Gamma(s)}\,\int_{0}^{1}\,\frac{f(y)-f(1/z)}{1-yz}\,\left|\log y\right|^{s-1}dy\label{eq:Borwein trick}\end{eqnarray}
 The second step is to find a sequence of polynomials $p_{n}(z)$
to be used for $f(z)$, such that the integral on the right hand side
can be evaluated as a simple, finite sum, while the left hand side
can be shown to bounded arbitrarily close to zero as $n\to\infty$.
The right hand side may be easily evaluated, by employing the following
theorem.

\begin{lem}
Given a polynomial $p_{n}(y)$ of degree $n$, it can be shown that
\begin{equation}
r_{n}(y)=\frac{p_{n}(y)-p_{n}(1/z)}{1-yz}\label{eq:r_n definition}\end{equation}
 is again a polynomial in $y$, of degree $n-1$, for any constant
$z$, provided that $z\ne0$. 
\end{lem}
\begin{proof}
This may be easily proved, term by term, by noting that $(x^{n}-a^{n})/(x-a)$,
with $a=1/z$, has the desired properties. 
\end{proof}
An explicit expression for $r_{n}(y)$ is needed. Write the polynomial
as \[
p_{n}(y)=\sum_{k=0}^{n}a_{k}y^{k}\]
 while for $r_{n}$ assume only a general series: \begin{equation}
r_{n}(y)=\sum_{k=0}^{\infty}c_{k}y^{k}\label{eq:r_n series}\end{equation}
 Setting $y=0$, one immediately obtains $c_{0}=a_{0}-p_{n}(1/z)$.
Higher coefficients are obtained by equating derivatives: \[
r_{n}^{(k)}(y)=\frac{1}{1-yz}\left[p_{n}^{(k)}(y)+kzr_{n}^{(k-1)}(y)\right]\]
 where $r_{n}^{(k)}(y)$ is the $k$'th derivative of $r_{n}$ with
respect to $y$. Setting $y=0$ in the above, one obtains the recurrence
relation $c_{k}=a_{k}+zc_{k-1}$ which is trivially solvable as \[
c_{k}=z^{k}\left[-p_{n}\left(\frac{1}{z}\right)+\sum_{j=0}^{k}\frac{a_{j}}{z^{j}}\right]\]
 From this, it is easily seen that $c_{n}=0$ and more precisely $c_{n+m}=0$
for all $m\ge0$. This confirms the claim that $r_{n}(y)$ is a polynomial
of degree $n-1$ in $y$. 

This is immediately employed in equation \ref{eq:Borwein trick} to
obtain \begin{equation}
\xi(s,z)=\mbox{Li}_{s}\left(z\right)+\frac{z}{p_{n}(1/z)}\sum_{k=0}^{n-1}\frac{c_{k}}{(k+1)^{s}}\label{eq:Polylog general estimate}\end{equation}
 To obtain a good approximation for $\mbox{Li}_{s}\left(z\right)$,
one needs to find a polynomial sequence such that $\xi$ goes to zero
as $n\to\infty$ for the domain $(s,z)$ of interest. That is, one
seeks to make \begin{equation}
\xi(s,z)=\frac{1}{p_{n}(1/z)}\,\frac{z}{\Gamma(s)}\,\int_{0}^{1}\,\frac{p_{n}(y)\left|\log y\right|^{s-1}}{1-yz}\, dy\label{eq:xi integral}\end{equation}
 as small as possible. One possible sequence, based on the Bernoulli
process or Gaussian distribution, is explored in a subsequent section.

\section{Convergence Acceleration of Oscillatory Series }

For the case of $z=-1$ and $p_{n}(y)=y^{n}$, the technique developed
above corresponds to Euler's method for the convergence acceleration
of an alternating series\cite[eqn 3.6.27]{A&S}; a particularly efficient
algorithm for Euler's method is given by van Wijngaarden\cite{Pre88}.
Polynomials that provide much more powerful series acceleration are
discussed by Cohen \emph{etal}\cite{Coh00}, but again in the context
of $z=-1$; these include the Tchebysheff orthogonal polynomials considered
by Borwein, which do a rather good job of minimizing the integrand
of equation \ref{eq:xi integral} when $z=-1$. In a certain sense,
an alternating series can be considered to be a series from which
an explicit factor of $z^{n}$ has been factored out, for the very
special case of $z=-1$. Thus, one is lead to ask if one can find
series acceleration methods for more general oscillatory series, where
the oscillatory component can be thought of as going as $z^{n}$ for
some root of unity $z$. 

The Euler method applied to the Hurwitz zeta function essentially
leads to the Hasse series \ref{eq:Hasse series}. This may be seen
by noting that the inner sum is just a forward difference:\[
\sum_{k=0}^{n}(-1)^{k}\left(\begin{array}{c}
n\\
k\end{array}\right)(q+k)^{1-s}=\triangle^{n}q^{1-s}\]
 where $\triangle$ is the forward difference operator. The Hasse
series has improved convergence in that the result is globally convergent,
whereas the traditional series representation of equation \ref{eq:Hurwitz zeta}
is not. However, a quick numerical experiment will show that the Hasse
series does not converge rapidly, especially in the critical strip.

The poor convergence of these traditional techniques motivates the
development of this paper. Strictly speaking, the slow convergence
can be attributed to the fact that the naive series representations
\ref{eq:Hurwitz zeta} or \ref{eq:polylog definition} are not strictly
alternating series, as would be required to justify the application
of the Euler method. Instead, the series are a superposition of oscillatory
terms; for the polylogarithm, an explicit oscillation coming from
$z^{n}=\left|z\right|^{n}e^{in\arg z}$ and a slower oscillation coming
from $n^{-s}=n^{-\sigma}e^{-i\tau\log n}$. Thus, what one needs is
a theory of sequence acceleration not for alternating series, but
for oscillatory series. The polylogarithm approximation represents
the \emph{ad-hoc} development of a special case for such a general
theory.

A loose sketch of a possible general theory is presented below. Suppose
that $c_{n}$ is some more or less arbitrary sequence of numbers,
oscillatory in $n$, and that one wishes to compute the sum \[
S=\sum_{n=0}^{\infty}c_{n}\]
 but that one wants to apply convergence acceleration techniques to
its evaluation. If the leading component of oscillation is $c_{n}\sim\cos n\theta$,
then one makes the Ansatz that there must exist some dual series $s_{n}$
oscillating as $s_{n}\sim\sin n\theta$ to give $e_{n}=c_{n}+is_{n}\sim e^{in\theta}$.
One may then consider a sum of $a_{n}=e^{in(\pi-\theta)}e_{n}$ and
consider the  summation of an alternating series \[
\sum_{n=0}^{\infty}\left(-1\right)^{n}a_{n}\]
 By construction, the $a_{n}$ are now presumed to be relatively tame,
so that series acceleration techniques can be applied. This trick
may be applied to any oscillatory sequence; one need not know a precise
or \emph{a priori} value for $\theta$; one need only guess at a reasonable
value, based on the data at hand. 

Alternately, one may consider replacing the forward difference operator
$\triangle b_{n}=b_{n+1}-b_{n}$ by $\triangle_{q}b_{n}=b_{n+1}-qb_{n}$.
Straightforward manipulations lead to \[
\sum_{m=0}^{\infty}\left(\frac{z}{1-qz}\right)^{m+1}\triangle_{q}^{m}b_{0}=\sum_{n=0}^{\infty}z^{n}b_{n}\]
 with the traditional Euler's series regained by taking $q=1$ and
$z=-1$. Assuming that the oscillatory part $z$ is known, or can
be approximately guessed at from the period of the oscillation, then
choosing $q=-1/z$ should make the left hand side an accelerated series
for the right hand side. However, there are more powerful techniques.

Suppose one is interested in the sum\begin{equation}
S(z)=\sum_{n=0}^{\infty}z^{n}b_{n}\label{eq: geo series}\end{equation}
 In general, one is interested in general sequences $b_{n}$. Suppose,
however, that the $b_{n}$ are expressible as moments, so that \[
b_{n}=\int_{0}^{1}y^{n}g(y)dy\]
 for some function $g(y)$. Since $g(y)$ is still fairly general
(subject to constraints discussed below), this assumption does not
overly restrict the $b_{n}$. If the $b_{n}$ are expressible as moments,
then one has\[
S(z)=\int_{0}^{1}\frac{g(y)}{1-yz}\, dy\]
 Applying the Borwein trick, one obtains \[
S(z)=S_{n}(z)+\xi(z)\]
 with \[
S_{n}(z)=\frac{1}{p_{n}\left(\frac{1}{z}\right)}\int_{0}^{1}\frac{p_{n}(y)-p_{n}\left(\frac{1}{z}\right)}{1-yz}g(y)\, dy\]
 and \[
\xi_{n}(z)=\frac{1}{p_{n}\left(\frac{1}{z}\right)}\int_{0}^{1}\frac{p_{n}(y)g(y)}{1-yz}\, dy\]
 The goal is to show that the $S_{n}(z)$ are easily evaluated, while
also showing that the $\xi_{n}(z)$ are arbitrarily small. The first
is easy enough: using equations \ref{eq:r_n definition} and \ref{eq:r_n series}
one has \[
S_{n}(z)=\frac{1}{p_{n}\left(\frac{1}{z}\right)}\sum_{k=0}^{n-1}c_{k}b_{k}\]
 With a suitable choice of $p_{n}(z)$, the coefficients $c_{k}$
are presumably easy to evaluate. To show that the $S_{n}(z)$ is really
a series acceleration for the geometric series \ref{eq: geo series},
one must show that $\left|\xi_{n}(z)\right|$ is bounded. If $g(y)\ge0$,
then \[
\left|\xi_{n}(z)\right|\le\left|\frac{p_{n}\left(y_{0}\right)}{p_{n}\left(\frac{1}{z}\right)}\right|\left|S(z)\right|\]
 where $y_{0}$ is the point on the unit interval where $p_{n}(y)$
assumes its maximum. Provided that one can find a polynomial sequence
such that \[
\left|\frac{p_{n}\left(y_{0}\right)}{p_{n}\left(\frac{1}{z}\right)}\right|\sim A^{-n}\]
 for some number $A>\left|1/z\right|$, then one has that the series
$S_{n}(z)$ converges to $S(z)$ more rapidly than the simple partial
sums $\sum_{k=0}^{n}z^{k}b_{k}$ of equation \ref{eq: geo series}.

Note that the above is a proof of a series acceleration method for
a more-or-less general sequence of $b_{n}$, and is in no way particular
to the polylogarithm or the Hurwitz zeta. The only ingredients of
the proof were that the $b_{n}$ are {}``well-behaved'' -- in this
case, being expressible as moments of a rather general function $g(y)$.
For the proof to hold, $g(y)$ needs to be (mostly) positive, and
integrable; it need not be analytic, differentiable or even continuous.
It needs to be {}``mostly'' positive only so that the bound on the
integrand, $p_{n}(y)g(y)\le\left|p_{n}\left(y_{0}\right)\right|g(y)$
can be established. Thus, one might hope that the acceleration method
might work, even for those general cases where $g(y)$ is not explicitly
known, but the $b_{n}$ are somehow {}``reasonable''. 

For the special case of $z=-1$, \emph{i.e.} for the case of an alternating
series, Cohen \emph{et al.}\cite{Coh00} suggest some remarkably strongly
converging polynomials, with values of $A$ from 5.8 to as much as
17.9. The generalization of those results to arbitrary $z$ is not
immediately apparent, but is surely possible. 

There is another, possibly fruitful direction one might take. This
hinges on noting that the integrands above take the general form of
a Padé approximation, and one proceeds from there.

\section{The Gaussian Distribution }

Returning to the polylogarithm, the task at hand is to find a polynomial
sequence that is suitable for establishing a tight bound on the error
term $\xi$. A good bound minimizes the integrand, while also maximizing
the value of $\left|p_{n}(1/z)\right|$. In order to suppress the
logarithmic branch point at $y=0$ in the integrand of \ref{eq:xi integral},
this section will consider the polynomial sequence $p_{n}(y)=y^{a}(1-y)^{n-a}$,
a sequence which goes over to the Gaussian distribution for large
$n$. 

Consider first the polynomial sequence $p_{2n}(y)=y^{n}(1-y)^{n}$.
It becomes sharply peaked at $y=1/2$, and subtends a diminishingly
small area. For large values of $n$, this polynomial sequence approximates
a Gaussian:\[
p_{2n}(y)=\frac{1}{4^{n}}e^{-4n(y-1/2)^{2}}\left[1+2\left(y-\frac{1}{2}\right)^{2}+\mathcal{O}\left(\left(y-\frac{1}{2}\right)^{4}\right)\right]\]
 Using this in the integral \ref{eq:xi integral}, and assuming that
the real part of $z$ is not positive, it is not hard to deduce the
very crude estimate\[
\xi(s,z)\sim\left(\frac{z^{2}}{4(z-1)}\right)^{n}\frac{z}{\Gamma(s)}\]
 which confirms that $\xi(s,z)$ does indeed get suitably small for
a certain region in the complex $z$-plane. However, in order for
equation \ref{eq:Polylog general estimate} to be useful computationally,
an upper bound on the value of $\xi$ needs to be given, as a function
of $z$ and $n$. This bound is derived in the next section. 

The polynomial coefficients are easily obtained, and are \begin{eqnarray*}
a_{0} & = & a_{1}=\cdots=a_{n-1}=0\\
a_{n+k} & = & \left(-1\right)^{k}\left(\begin{array}{c}
n\\
k\end{array}\right)\quad\mbox{ for }0\le k\le n\end{eqnarray*}
 This leads to \[
\mbox{Li}_{s}(z)=-\frac{z^{2n+1}}{\left(z-1\right)^{n}}\sum_{k=0}^{2n-1}\frac{c_{k}}{\left(k+1\right)^{s}}+\xi(s,z)\]
where the $c_{k}$ are given by \[
c_{k}=z^{k}\left[-\left(\frac{z-1}{z^{2}}\right)^{n}+\frac{1}{z^{n}}\sum_{j=0}^{k-n}\left(-1\right)^{j}\left(\begin{array}{c}
n\\
j\end{array}\right)\frac{1}{z^{j}}\right]\]
 where the summation above is to be understood to be zero when $k<n$.
The above summations can be re-organized into the more suggestive
form \begin{equation}
\mbox{Li}_{s}(z)=\sum_{k=1}^{n}\frac{z^{k}}{k^{s}}\;+\frac{1}{\left(1-z\right)^{n}}\;\sum_{k=n+1}^{2n}\frac{z^{k}}{k^{s}}\;\sum_{j=0}^{2n-k}\left(-z\right)^{j}\left(\begin{array}{c}
n\\
j\end{array}\right)+\xi(s,z)\label{eq:polylog specific estimate}\end{equation}
 The error term $\xi(s,z)$ is negligible for only a very specific
area of the complex $z$-plane. The region where $\xi$ may be ignored
as a small error is derived in the next section. Substituting $z=-1$
in the above formulas agrees with expressions previously given by
Borwein\cite{Bor95}. 

The above formula has been implemented using the Gnu Multiple Precision
arithmetic library (GMP)\cite{GMP}, and has been numerically validated
for correctness in several different ways. The source code, under
the license terms of the LGPL, may be obtained by contacting the author.

\section{Bound on the Error Term}

In order for equation \ref{eq:polylog specific estimate} to be useful
computationally, an upper bound on the value of $\xi$ as a function
of $n$ needs to be given. To compute the polylogarithm to some desired
precision, one infers a suitable value of $n$ based on this bound.
However, to achieve this desired precision, one must not only choose
$n$ large enough, but one must also maintain a somewhat larger number
of digits in the intermediate terms, as the appearance of the binomial
coefficient in the equation \ref{eq:polylog specific estimate} implies
that intermediate terms can become quite large, even while mostly
cancelling. This section derives an upper bound on the size of $\xi$,
and briefly discusses the issue of the required precision in intermediate
terms. 

The general behavior of the integrand appearing in equation \ref{eq:xi integral}
depends on the whether $\Re s\ge1$ or not; estimates are presented
for these two cases. Writing $s=\sigma+i\tau$ for real $\sigma$
and $\tau$, and assuming $\sigma\ge1$ and choosing $n$ so that
$\sigma\le n$, one has\begin{eqnarray*}
\left|p_{2n}(y)\left|\log y\right|^{s-1}\right| & \le & p_{2n}(y)\left|\log y\right|^{\sigma-1}\\
 & = & \left|y(1-y)\log y\right|^{\sigma-1}\left(y(1-y)\right)^{n-\sigma+1}\end{eqnarray*}
 Each part of the right hand side is bounded by a Gaussian, and since
the product of Gaussians is a Gaussian, so is the entire expression.
From the Stirling approximation, one has the well-known identity \[
\left(y(1-y)\right)^{a}\le\frac{1}{4^{a}}\exp-4a\left(y-\frac{1}{2}\right)^{2}\]
 The other part is also bounded by a Gaussian \[
\left|y(1-y)\log y\right|^{a}\le\left|y_{0}\left(1-y_{0}\right)\log y_{0}\right|^{a}\exp\left(-\frac{a\left(y-y_{0}\right)^{2}}{2y_{0}\left(1-2y_{0}\right)}\right)\]
 which is centered on the maximum of $\left|y_{0}\left(1-y_{0}\right)\log y_{0}\right|$,
that is, at $y_{0}=0.23561058253\ldots$, where $y_{0}$ is the solution
to \[
0=\frac{d}{dy}\; y(1-y)\left|\log y\right|\]
 Curiously, the numerical value of the factor in the Gaussian is $1/2y_{0}(1-2y_{0})=4.013295587\ldots$

To evaluate the integral in \ref{eq:xi integral}, one also needs
a bound on the denominator. This is furnished by \[
\left|\frac{1}{1-yz}\right|\le C(z)=\begin{cases}
1 & \mbox{ if }\Re z\le0\\
\frac{1}{\left|1-z\right|} & \mbox{ if }0<\Re z\;\mbox{ and }\left|z\right|<1\\
\left|\frac{z}{\Im z}\right| & \mbox{ if }0<\Re z\;\mbox{ and }\left|z\right|>1\end{cases}\]
 Thus, the integrand is bounded by a Gaussian. Multiplying the two
Gaussians, completing the square, and evaluating the resulting integral
is a bit tedious. The result is that \[
\left|\int_{0}^{1}\,\frac{p_{2n}(y)\left|\log y\right|^{s-1}}{1-yz}\, dy\right|\le C(z)\frac{\left|4y_{0}\left(1-y_{0}\right)\log y_{0}\right|^{\sigma-1}}{4^{n}}\, G(s)\]
 The constant may be taken as $\left|4y_{0}\left(1-y_{0}\right)\log y_{0}\right|=1.041381965\ldots$.
The factor $G(s)$ is bounded by \begin{eqnarray*}
G(s) & \le & \sqrt{\frac{\pi}{4n}}\exp-\left(1-2y_{0}\right)^{2}\left(\sigma-1\right)\frac{8y_{0}\left(1-2y_{0}\right)}{8y_{0}\left(1-2y_{0}\right)+\frac{\sigma-1}{n-\sigma+1}}\\
 & \le & \sqrt{\frac{\pi}{4n}}\exp-\left(1-2y_{0}\right)^{2}\frac{\sigma-1}{\sigma}\\
 & \le & 1\end{eqnarray*}
 when $1\le\sigma$. Useful for estimation is that $8y_{0}(1-2y_{0})=0.996687115\ldots$.
Combining these results, one obtains \begin{equation}
\left|\xi(s,z)\right|\le\frac{1}{4^{n}}\left|\frac{z^{2}}{z-1}\right|^{n}\left|\frac{z}{\Gamma(s)}\right|C(z)G(s)\left(1.0414\ldots\right)^{\sigma-1}\label{eq:error-bound-positive}\end{equation}
 for the case where $1\le\sigma$ and $n$ chosen so that $\sigma\le n$.

The case where $\sigma\le1$ may be evaluated as follows. One writes
\begin{eqnarray*}
\left|p_{2n}(y)\left|\log y\right|^{s-1}\right| & \le & p_{2n}(y)\left|\log y\right|^{\sigma-1}\\
 & = & \left|\frac{1-y}{\log y}\right|^{1-\sigma}y^{n}(1-y)^{n+\sigma-1}\\
 & \le & y^{n}(1-y)^{n+\sigma-1}\end{eqnarray*}
 The integrand may now be recognized as the Beta function, so that\[
\left|\int_{0}^{1}\,\frac{p_{2n}(y)\left|\log y\right|^{s-1}}{1-yz}\, dy\right|\le C(z)\frac{\Gamma(n+1)\Gamma(n+\sigma)}{\Gamma(2n+\sigma+1)}\]
 In the region where $\left|\sigma\right|\ll n$, one may approximate
\[
\frac{\Gamma(n+1)\Gamma(n+\sigma)}{\Gamma(2n+\sigma+1)}\le\frac{2^{1-\sigma}}{4^{n}}\]
 Combining the various pieces, one obtains a result remarkably similar
to the previous one; namely, when $\sigma\le1$ but $n$ is such that
$\left|\sigma\right|\ll n$, one has \begin{equation}
\left|\xi(s,z)\right|\le\frac{1}{4^{n}}\left|\frac{z^{2}}{z-1}\right|^{n}\left|\frac{z}{\Gamma(s)}\right|2^{1-\sigma}C(z)\label{eq:error-bound-negative}\end{equation}

For evaluations on the critical line $\sigma=1/2$, one needs a good
estimate for $\left|\Gamma(s)\right|^{-1}$, which can become quite
large as the imaginary part of $s$ increases. A useful estimate for
$\left|\Gamma(s)\right|^{-1}$ is given by Borwein \cite{Bor00-A,Bor95},
for the case where $\sigma\ge-m+1/2$: \begin{eqnarray*}
\frac{1}{\left|\Gamma(s)\right|} & = & \frac{1}{\left|\Gamma(\sigma)\right|}\sqrt{\prod_{k=0}^{\infty}\left(1+\frac{\tau^{2}}{\left(\sigma+k\right)^{2}}\right)}\\
 & \le & \frac{1}{\left|\Gamma(\sigma)\right|}\sqrt{\frac{\sinh\tau\pi}{\tau\pi}\,\prod_{k=0}^{m}\left(1+\frac{\tau^{2}}{\left(\frac{1}{2}+k\right)^{2}}\right)}\end{eqnarray*}

A remarkable side-effect of the estimation is that, for $s$ equal
to negative integers, one has that $\Gamma(s)$ is infinite, thus
implying that the error term is zero. This somewhat surprising result
corresponds to the fact that $\mbox{Li}_{-n}(z)$ has an exact expression
as a ratio of two polynomials, the denominator being of degree $n+1$.
Indeed, setting $s=0$ and so $n=1$ into equation \ref{eq:polylog specific estimate}
gives the exact result \[
\mbox{Li}_{0}(z)=\frac{z}{1-z}\]
 while for $s=-1$, one must use $n=2$, to obtain the exact result\[
\mbox{Li}_{-1}(z)=\frac{z}{\left(1-z\right)^{2}}\]
 The generic form for the polylogarithm at the negative integers is
\[
\mbox{Li}_{-n}(z)=\left(-1\right)^{n}\,\sum_{k=0}^{n}\left(-1\right)^{k}\frac{k!}{\left(1-z\right)^{k+1}}\left\{ \begin{array}{c}
n+1\\
k+1\end{array}\right\} \]
 where $\left\{ \begin{array}{c}
n\\
k\end{array}\right\} $ denotes the Stirling numbers of the second kind.

In general, the formula \ref{eq:polylog specific estimate} seems
to have the curious property of always evaluating to exactly the same
rational function whenever $s$ is a non-positive integer, and $-s<n$.
That is, for fixed non-positive integer $s$, the sum is independent
of $n$, provided that $n$ is large enough. Thus, this polynomial
approximation has the pleasing result of giving exact answers in the
case when an exact polynomial answer is possible.

For the general case, one may conclude that the estimate \ref{eq:polylog specific estimate}
can be effectively applied whenever\[
\left|\frac{z^{2}}{z-1}\right|<4\]
 To obtain a value of the polylogarithm to within some given numerical
precision, one must invert the formulas \ref{eq:error-bound-positive}
or \ref{eq:error-bound-negative}, solving for the value of $n$ which
is to be used in equation \ref{eq:polylog specific estimate}. To
obtain a fixed number of digits of precision, one must carry out intermediate
calculations with at least that many digits of precision. In fact,
one must have more. The appearance of the binomial coefficient in
equation \ref{eq:polylog specific estimate} is the problem. Since
$2^{n}>\left(\begin{array}{c}
n\\
j\end{array}\right)$, one concludes that $n$ additional binary bits of precision are
needed with which to carry out the calculations. As a result, it becomes
difficult to implement the algorithm robustly using only double-precision
arithmetic. For values of $z$ near $z=-1$, the algorithm is well-enough
behaved; however, round-off errors significantly disturb the calculations
as $z$ approaches $+1$. This trouble with precision can be at least
partly alleviated by making use of the duplication formula, as discussed
in the next section.

\section{Euler-Maclaurin Summation}

An alternative algorithm for computing the Hurwitz zeta may be based
upon the Euler-Maclaurin series. This algorithm proves to be quite
rapid and efficient\cite{Pav07}. The Euler-Maclaurin series may be
written as\begin{eqnarray}
\sum_{k=0}^{\infty}f(k) & = & \sum_{k=0}^{N-1}f(k)+\frac{f(N)}{2}+\int_{N}^{\infty}f(x)dx\nonumber \\
 &  & -\sum_{k=1}^{p}\left.\frac{B_{2k}}{(2k)!}\,\frac{d^{2k+1}}{dx^{2k+1}}f(x)\right|_{x=N}+R\label{eq:Euler-Maclaurin}\end{eqnarray}
 This formula is simply applied to the function $f(x)=(x+q)^{-s}$
to obtain the Hurwitz zeta. The derivative is particularly easy to
evaluate: \[
\frac{d^{2k+1}}{dx^{2k+1}}\,\frac{1}{(x+q)^{s}}=-\frac{s(s+1)(s+2)\cdots(s+2k)}{(x+q)^{s+2k+1}}\]
 and the integral is trivial:\[
\int_{N}^{\infty}\frac{1}{(x+q)^{s}}\, dx=\frac{1}{(s-1)}\;\frac{1}{(N+q)^{s-1}}\]
 The error made by this expansion is embodied in the term $R$. It
is directly given by \[
R=\frac{2}{(2\pi)^{2p}}\int_{N}^{\infty}f^{(p+1)}(x)\frac{B_{p+1}\left(x-\left\lfloor x\right\rfloor \right)}{(p+1)!}dx\]
where the $B_{k}(x)$ are the Bernoulli polynomials. There are various
formulas in the literature which may be readily applied to bound the
size of this term\cite{A&S}. Thus, for example, one has \[
R\le\frac{2}{(2\pi)^{2p}}\frac{\left|s(s+1)\cdots(s+2p)\right|}{\Re s+2p}\,\frac{1}{(N+q)^{\Re s+2p}}\]
 Thus, the evaluation of the sum \ref{eq:Euler-Maclaurin} merely
requires a suitable choice of $N$ and $p$ to be used. As it happens,
the situation is even simpler than this. The sum \ref{eq:Euler-Maclaurin}
is an asymptotic series, and it is well-known that the best approximation
for an asymptotic series occurs when one stops the summation at the
smallest term in the series. Thus, it is sufficient to choose a value
of $N$, while $p$ is found dynamically as the algorithm runs. 

There remains the question of what value of $N$ to use. This remains
an open problem, which will not be explored here. However, it seems
that an adequate heuristic for most common cases is to choose $N=D/2+10$
where $D$ is the desired number of decimal digits of precision. The
performance of this algorithm is compared to that of the Borwein algorithm
in the next section.

XXX to do: characterize the region of convergence for this algo.

\section{Performance}

\begin{figure}

\caption{\label{cap:Polylogarithm-cold-cache}Polylogarithm cold cache performance}

\includegraphics[%
  width=1\textwidth]{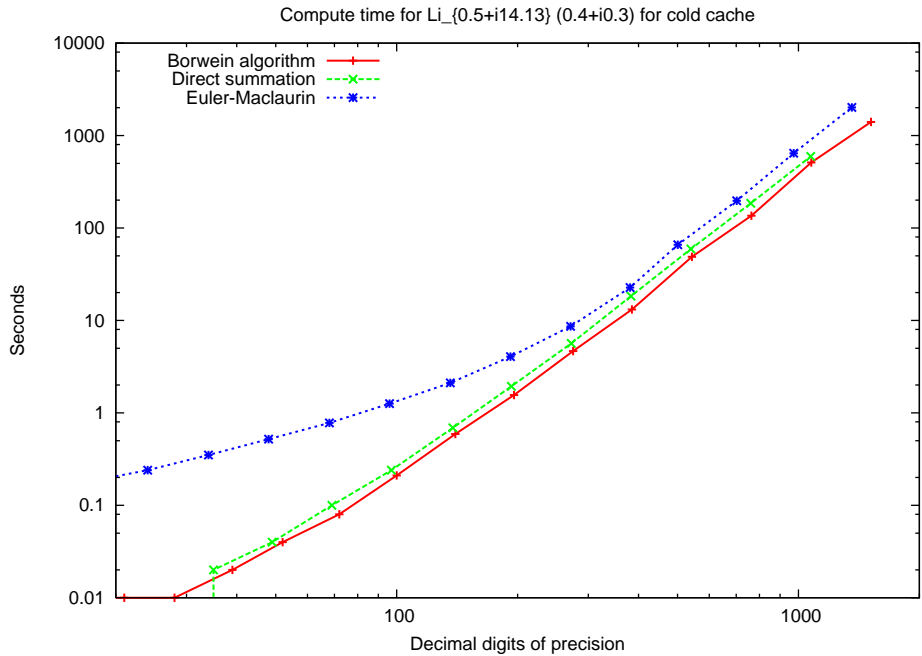}

This figure shows the time, in seconds, needed to compute the polylogarithm
to the indicated number of decimal places of precision. The values
chosen for evaluation are $s=0.5+i14.134725$, which is very near
to the first of the non-trivial Riemann zeroes, and $z=0.4+i0.3$,
which is relatively near to the origin. For such a small value of
$z$ ($\left|z\right|=0.5$), one might hope that direct summation
might compete favorably against the other two algorithms; and in particular
against the Borwein algorithm, as the point is closer to the troublesome
branch point at $z=+1$ than algorithmically optimum $z=-1$. Nevertheless,
$\left|z^{2}/(z-1)\right|\approx0.373$ for this point, and the Borwein
algorithm wins. The Euler-Maclaurin algorithm is relatively disadvantaged,
as it needs to be evaluated twice, and a value for $\Gamma(s)$ must
be computed as well. 

The primary contribution to the calculations is the evaluation of
$n^{-s}$ in the summations. If these are pre-computed, the Borwein
algorithm still wins, as shown in the next figure. All calculations
were performed on a contemporary desktop computer, using the GMP library;
the algorithms were implemented in the C programming language. Note
that the actual performance depends on the underlying implementation
of log, exp, sin, atan, sqrt, gamma and the like; the actual implementation
used here is fairly naive and untuned.
\end{figure}

\begin{figure}

\caption{\label{cap:Polylogarithm-warm-cache}Polylogarithm warm cache}

\includegraphics[%
  width=1\textwidth]{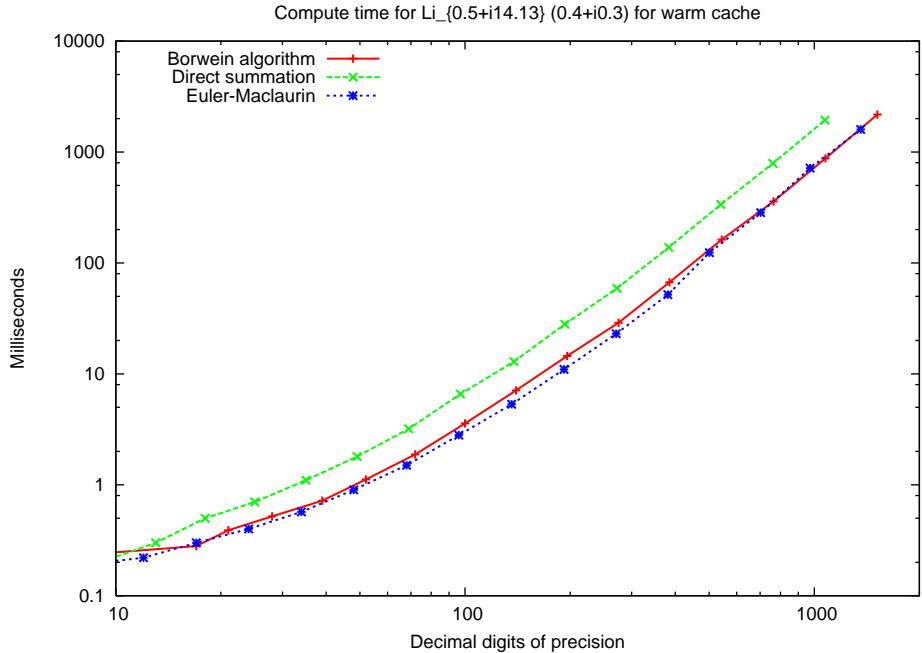}

This graph compares the compute times for the Borwein algorithm, the
direct summation of the polylogarithm, and the Euler-Maclaurin series,
for a {}``warm cache''. That is, an array of values of $n^{-s}$
were pre-computed, prior to beginning the timing measurements. These
pre-computed values were used in the Borwein and direct sums. For
the Euler-Maclaurin sum, the value of $\Gamma(s)/(2\pi)^{s}$ is pre-computed
and cached as well. Despite the use of such cached, pre-computed values,
the Borwein algorithm still wins over direct summation, and appears
to be tied with the Euler-Maclaurin sum. This graph uses the same
$s$ and $z$ values as in the previous graph. 
\end{figure}

The performance of both the Euler-Maclaurin summation and the Borwein
algorithm appears to be very good. As a general rule, both have better
performance than the direct summation formula \ref{eq:polylog definition},
although this depends on whether the polylogarithm or the Hurwitz
zeta is the desired quantity. It also depends on whether repeated
evaluations are being done for fixed $s$ but varying $q$ or $z$,
as holding $s$ fixed allows the results of intermediate calculations
to be cached. The figures \ref{cap:Polylogarithm-cold-cache} and
\ref{cap:Polylogarithm-warm-cache} compare the performance of actual
implementations.

\begin{figure}

\caption{Hurwitz zeta cold cache performance}

\includegraphics[%
  width=1\textwidth]{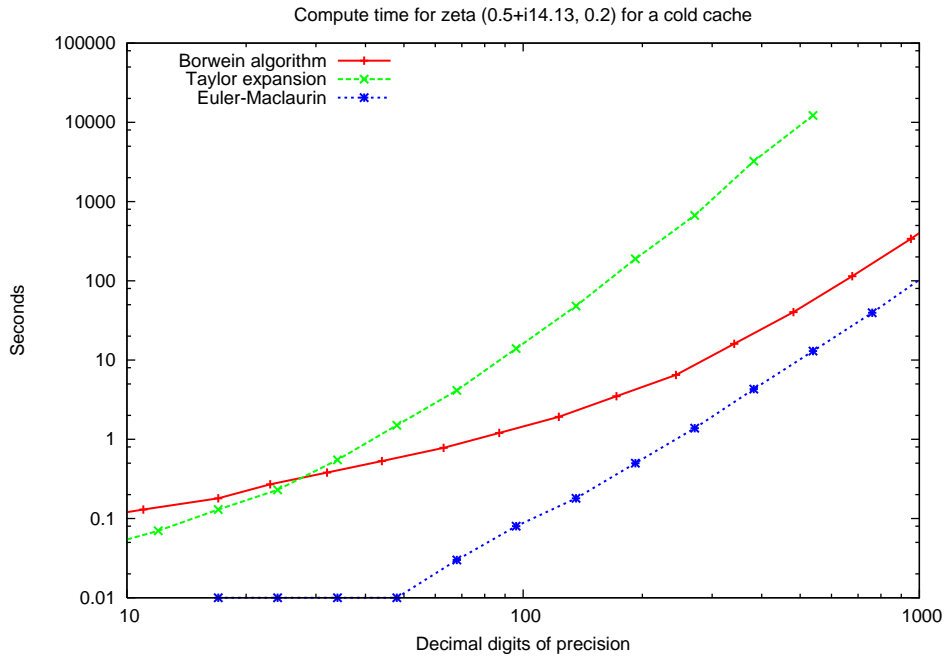}

This figure compares the performance of the Borwein algorithm, the
evaluation of the Taylor expansion, and Euler-Maclaurin summation,
for the Hurwitz zeta. The vales computed are for $s=0.5+i14.13$ and
$q=0.2$. For such a small value of $q$, one might have hoped that
the Taylor expansion might converge quickly. This appears to not be
the case, as the binomial coefficients can grow to be quite large,
and thus force many terms to be summed. 

The figure is for a {}``cold cache'', assuming that no constants
have been pre-computed. The lines are not parallel, with the Taylor
expansion getting progressively worse. Both the Borwein and the Euler-Maclaurin
algorithms are orders of magnitude faster than the Taylor's series
for high precisions; the Euler-Maclaurin algorithm appears to be an
easy winner for all precisions.
\end{figure}

\begin{figure}

\caption{Hurwitz zeta warm cache performance}

\includegraphics[%
  width=1\textwidth]{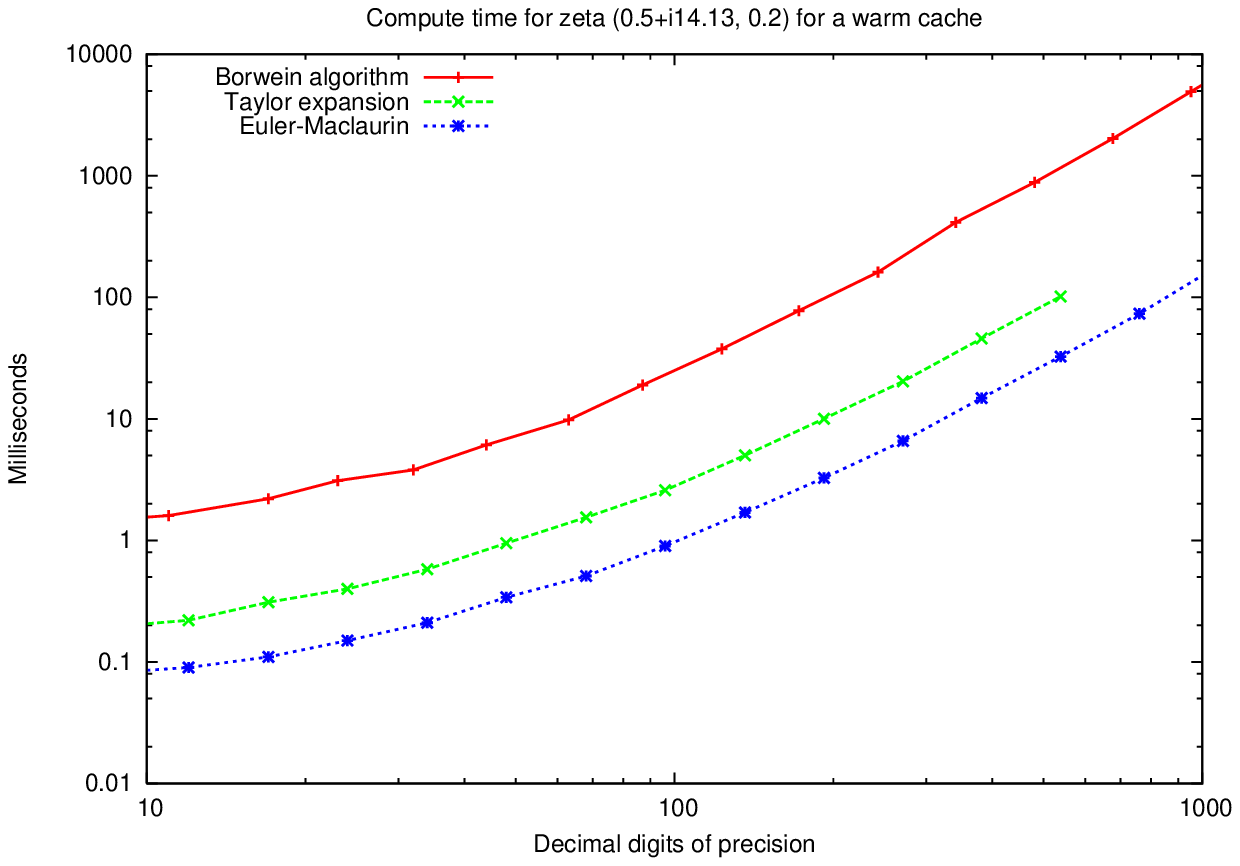}

This figure compares the evaluation times for Euler-Maclaurin summation,
the Borwein algorithm, and the Hurwitz zeta Taylor's expansion, for
the {}``warm cache'' scenario. In this case, it is recognized that
if $s$ is assumed to be held fixed, then the values of $n^{-s}$
appearing in the Borwein algorithm may be pre-computed. Similarly,
the values of $\zeta(s+n)$ and the binomial coefficients $\left(\begin{array}{c}
s+n-1\\
n\end{array}\right)$ appearing in the Taylor's series and the Euler-Maclaurin formula
may be pre-computed.

For this value of $q$, the Borwein algorithm requires four evaluations
of the polylogarithm; values of $q$ closer to 0.5 would require only
two. The Taylor's series evaluation appears to be about eight times
faster than these four evaluations (or four times faster than the
minimum of two evaluations).

Note that the time axis here is in milliseconds, not seconds. Comparing
to the previous graph, it is clear that performing the pre-computations
can be terribly expensive. The first evaluation of the function can
take 100 or 1000 times longer than subsequent evaluations at the same
value of $s$.
\end{figure}

The first figure shows {}``cold cache'' performance, measured in
seconds, as compared to the number of desired decimal places of precision.
The figure is termed {}``cold cache'', in recognition of the fact
that some constants may be pre-computed. For example, if $s$ is held
fixed, while $z$ is varied, then the values of $n^{-s}$ appearing
in both the direct sum and the Borwein algorithm may be computed once,
and then re-used for subsequent calculations. As the figures indicate,
computing $n^{-s}$ can be very expensive for general complex-valued
$s$, and so the caching strategy offers a big gain when $s$ is held
constant. 

As the {}``cold cache'' figure demonstrates, the Borwein algorithm
is faster than direct summation. The problem with direct summation
is that it requires more values of $n^{-s}$ to be computed to achieve
comparable precision. 

As the {}``warm cache'' figure shows, the Borwein algorithm will
still win against direct summation, even when the values of the $n^{-s}$
are pre-computed. Even when these are pre-computed, and can be pulled
from the cache, direct summation still requires more operations.

A different but equally dramatic set of results obtain, when one compares
the performance of the Borwein algorithm applied to the Hurwitz zeta
function, as compared to the use of the Taylor's expansion \ref{eq:Hurwitz taylor}
for the Hurwitz zeta. In this case, as the figures show, the Taylor's
series outperforms the Borwein algorithm by a constant factor of three,
when both algorithms use pre-computed constants. However, the cost
of pre-computing these constants skyrockets for the Taylor's series,
making it unattractive as the number of required decimal places increases.

Although the Borwein algorithm seems to show a slight advantage over
the Euler-Maclaurin series when it is used to evaluate the polylogarithm,
much of that advantage disappears when evaluating near the branch
point $z=1$. The Borwein algorithm breaks down near the branch point,
and requires the (possibly recursive) application of the duplication
formula (given in the next section) to obtains points sufficiently
far away from $z=1$. By contrast, the Euler-Maclaurin series seems
to happily converge in this area, and so no extra steps are required. 

In conclusion, it appears that the Euler-Maclaurin formula provides
the best algorithm for evaluating the Hurwitz zeta, and can often
tie and sometimes outperform the Borwein algorithm for the polylogarithm.

\section{Duplication formula}

The region of applicability of the algorithm may be extended by making
use of the duplication formula for the polylogarithm or periodic zeta.
The duplication formula, or more generally, the multiplication theorem,
is an extension of the well-known Legendre duplication formula for
the Gamma function\cite[(6.1.18),(6.3.8)(23.1.10)]{A&S}, into the
domain of the polylogarithm.

Thus, for example, the formula \ref{eq:polylog specific estimate}
together with the error bounds \ref{eq:error-bound-positive},\ref{eq:error-bound-negative}
allow $F(q;s)=\mbox{Li}_{s}\left(e^{2\pi iq}\right)$ to be computed
for real $q$ in region $1/4\le q\le3/4$. To obtain values on the
region $0<q<1/4$, one applies the duplication formula \[
F\left(q;s\right)=2^{1-s}F\left(2q;s\right)-F\left(q+\frac{1}{2};s\right)\]
 recursively until one obtains values of $q\ge1/4$. Rearranging terms,
one obtains a similar formula for iterating values of $q>3/4$ until
they are less than or equal to $3/4$. Thus, as $q$ approaches 0
or 1, the algorithm requires more time, but only logarithmically so,
as $-\log_{2}q$ or $-\log_{2}(1-q)$ as the case may be.

The duplication formula follows from a general $p$-adic identity
\[
\sum_{m=0}^{p-1}F\left(\frac{m+q}{p};s\right)=p^{1-s}F\left(q;s\right)\]
 which is valid for any positive integer $p$. One of many ways of
obtaining this formula is by noting that the function $F\left(s;q\right)$
is an eigenvalue of the $p$-adic Bernoulli transfer operator associated
with eigenvalue $p^{1-s}$\cite{Ve-B04}.

The equivalent formulas\cite[Sections 7.3.1, 7.12.1]{Lew81} for the
polylogarithm are \[
\mbox{Li}_{s}(z)+\mbox{Li}_{s}(-z)=2^{1-s}\mbox{Li}_{s}\left(z^{2}\right)\]
 whereas the $p$-adic version has the appearance of a Gauss sum:\[
\sum_{m=0}^{p-1}\mbox{Li}_{s}\left(ze^{2\pi im/p}\right)=p^{1-s}\mbox{Li}_{s}\left(z^{p}\right)\]

The application of the duplication formula, together with the inversion
relation \ref{eq:poly-hur relation} can be used to extend the evaluation
of the polylogarithm to the entire complex plane. For numerical work,
both formulas must be applied, one alone is not enough. Consider first
the application of the duplication formula only. It may be used to
take points that are near to $z=1$, and map them further away from
$z=1$, into the convergent region for the Borwein polynomial. Repeated
application allows arbitrarily close approach to $z=1$ from the left-hand
side. The resulting region of convergence is kidney-shaped, with the
cusp of the kidney at $z=1$, and the kidney containing the unit disk
$\left|z\right|\le1$. The precise shape of the kidney depends on
the number of terms one wishes to use in the polynomial approximation.
The shape that can be achieved while still maintaining good running
time is shown in figure \ref{cap:Polylogarithm-phase-plot}. However,
as can be seen, this strategy barely penetrates the $\Re z>1$ region.

To extend the algorithm to the entire complex $z$-plane, one must
make use of the Jonquière's inversion formula\cite[pp 27-32]{Bate53}

\begin{equation}
e^{-i\pi s/2}\mbox{Li}_{s}\left(z\right)+e^{i\pi s/2}\mbox{Li}_{s}\left(\frac{1}{z}\right)=\frac{\left(2\pi\right)^{s}}{\Gamma(s)}\zeta\left(1-s,\frac{\log z}{2\pi i}\right)\label{eq:Inversion formula}\end{equation}
 together with some sort of independent means of evaluating the Hurwitz
zeta function. The Taylor expansion \ref{eq:Hurwitz taylor} is particularly
well-suited, as it is rapidly convergent in the vicinity of $z=1$.
Specifically, it is convergent when $\left|\log z\right|<2\pi$, although
the region of acceptable numerical convergence is smaller, roughly
$\left|\log z\right|<\pi$. Either way, this encompasses a rather
large region in the vicinity of $z=1$, which is exactly the region
where the Borwein algorithm has trouble. The inversion formula is
used to move points $\left|z\right|>1$ from outside of the unit circle,
to the inside. Not all points inside the unit circle are directly
accessible to the Borwein algorithm; the duplication formula must
still be used to handle points in the disk interior that are near
$z=1$. Similarly, points with very large $z$, such as those for
which $\left|z\right|>e^{2\pi}$, require the use of the duplication
formula to be moved into the region of convergence for the Hurwitz
Taylor series. A typical view of the polylogarithm on the complex
$z$-plane is shown in figure \ref{cap:Polylog-whole-plane}.

As a practical computational matter, rather than using \ref{eq:Inversion formula},
it seems to be more convenient to use the formula\begin{equation}
\mbox{Li}_{1-s}\left(z\right)=\frac{\Gamma(s)}{\left(2\pi\right)^{s}}\left[e^{i\pi s/2}\zeta\left(s,\frac{\log z}{2\pi i}\right)+e^{-i\pi s/2}\zeta\left(s,1-\frac{\log z}{2\pi i}\right)\right]\label{eq:polylog as pure hurwitz}\end{equation}
 Although this is equivalent to the equation \ref{eq:Inversion formula}
immediately above, it correctly captures the appropriate branch of
the logarithm that should be used: equation \ref{eq:polylog as pure hurwitz}
works well for values of $s$ in both the upper and lower half-planes,
whereas equation \ref{eq:Inversion formula} is more difficult to
correctly evaluate in the lower half $s$-plane. In either case, one
must use the logarithm with an unusual branch cut, taking it to run
from $z=0$ to the right along the positive real axis, as opposed
to the usual cut taken for the logarithm. 

It is of some curiosity to note that the sheet that is forced takes
the form \[
\zeta\left(s,q\right)+e^{-i\pi s}\zeta\left(s,1-q\right)=\sum_{n=-\infty}^{+\infty}\frac{1}{(n+q)^{s}}\]
 which is also noted by Costin\cite[eqn 14]{Cos07} as a Mittag-Leffler
type decomposition. A general discussion of the monodromy is given
in a later section.

\section{Testing and Validation}

The correctness of a given numerical implementation can be validated
in a number of ways. For non-positive integer $s$, one has the exact
rational functions previously mentioned. For positive integer $n$,
one has the relationship \[
\mbox{Li}_{n}\left(e^{2\pi iq}\right)+\left(-1\right)^{n}\mbox{Li}_{n}\left(e^{-2\pi iq}\right)=-\frac{\left(2\pi i\right)^{n}}{n!}B_{n}(q)\]
 where the $B_{n}(x)$ are the Bernoulli polynomials. For $z=-1$,
one regains the Riemann zeta function $\zeta(s)$: \[
\mbox{Li}_{s}\left(-1\right)=\frac{1}{2^{1-s}-1}\,\zeta(s)\]
 For $\left|z\right|<1$, the defining series \ref{eq:polylog definition}
is explicitly convergent, and may be directly summed, thus offering
a fourth check of the correctness of an implementation. Finally, the
multiplication theorem can be used to check for consistency. Each
of these quantities can be computed by independent algorithms, and
thus be used to validate the correctness of a polylogarithm implementation.

\section{Branch Points and Monodromy}

The principal sheet of the polylogarithm has a branch point at $z=1$,
and by convention, a branch cut is placed along the positive real
$z$-axis, extending to the right. As one moves off of the principal
sheet, one discovers that there is another branch point at $z=0$.
The resulting monodromy group is generated by two elements, acting
on the covering space of the bouquet $S^{1}\vee S^{1}$ of homotopy
classes of loops in $\mathbb{C}\setminus\left\{ 0,1\right\} $ passing
around the branch points $z=0$ or $z=1$. The author is not aware
of any simple published discussion of the monodromy for the polylogarithm;
a much more abstract discussion is given in \cite{Cos07,Hain92,Hain91,Bloch91}.
Thus, this section provides a discussion; note that the correct manipulations
to move from one sheet to another can be somewhat treacherous and
confusing. The 

\begin{figure}

\caption{Polylogarithm Monodromy}

\includegraphics[%
  width=1\textwidth]{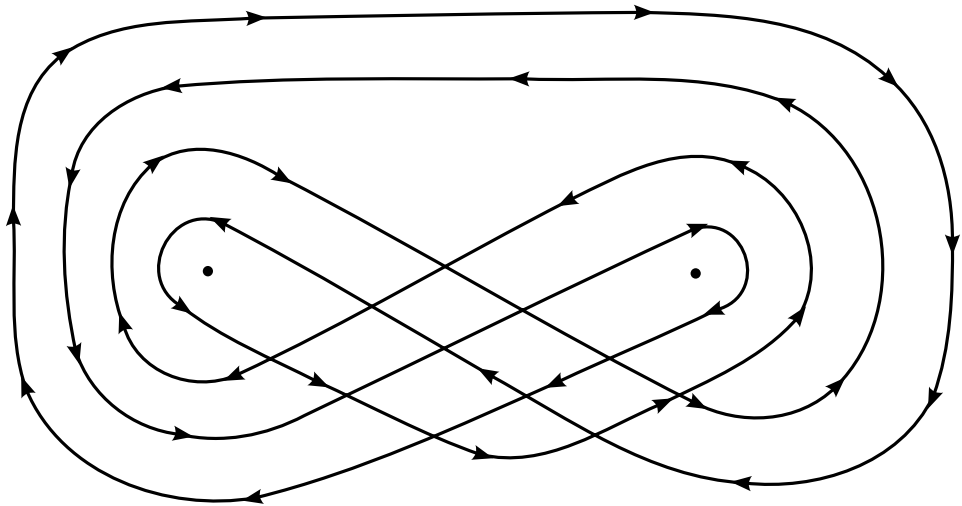}

A graph of the polylogarithm monodromy, given by equation \ref{eq:polylog monodromy}.
\end{figure}

The inversion formula \ref{eq:Inversion formula} suggests a way to
move around the branch point at $z=1$. Suppose one starts at the
point $z=x+i\epsilon$ with $x$ real, positive, and greater than
one, and $\epsilon$ some arbitrarily small positive real number;
thus $z$ is very near the branch cut of the principal sheet. One
wishes to compare this to the straddling value at $x-i\epsilon$.
Applying the inversion formula, one can bounce these two points inside
the unit circle, where there is no cut, and thus the polylog differs
by $\mathcal{O}(\epsilon)$. The difference across the cut is thus
\begin{eqnarray*}
\Delta & = & \mbox{Li}_{s}\left(x+i\epsilon\right)-\mbox{Li}_{s}\left(x-i\epsilon\right)\\
 & = & e^{i\pi s/2}\frac{\left(2\pi\right)^{s}}{\Gamma(s)}\left[\zeta\left(1-s,\frac{\ln x+i\epsilon}{2\pi i}\right)-\zeta\left(1-s,\frac{\ln x-i\epsilon}{2\pi i}\right)\right]+\mathcal{O}(\epsilon)\end{eqnarray*}
 Now, since $x>1$, a naive application of this formula yields $\Delta=0$
since $\log(x+i\epsilon)-\log(x-i\epsilon)=\mathcal{O}(\epsilon)$;
but clearly $\Delta$ cannot be zero. To resolve this situation, one
must make the Ansatz that the cut of the logarithm should be crossed.
This may be done by taking the cut of the logarithm to extend to the
right, instead of to the left, as it is by convention. One then gets
\begin{equation}
\log(x+i\epsilon)-\log(x-i\epsilon)=2\pi iN+\mathcal{O}(\epsilon)\label{eq:log-discontinuity}\end{equation}
 for some integer $N$. The difference across the cut, for positive
$N$, is then \[
\zeta\left(s,q\right)-\zeta\left(s,N+q\right)=\sum_{k=0}^{N-1}\frac{1}{\left(k+q\right)^{s}}\]
 with $q=\ln(x)/2\pi i$. By swinging the cut of the logarithm so
that it extends to the right, one obtains that the real part of $q$
is positive; the real part of $q$ runs between $0$ and $1$. As
a result, the value of $\left(k+q\right)^{s}$ is unambiguous; as
otherwise taking something to the power $s$ also begs the question
of which sheet the power must be taken on. That is, for general complex
$w$, the power function is multi-sheeted:\begin{equation}
w^{s}=e^{s\ln w}\to e^{s(\ln w+2\pi iM)}=e^{2\pi isM}w^{s}\label{eq:M-sheet}\end{equation}
 for some integer $M$ representing the sheet of the power function.
Since the real part of $q$ is taken as positive, one can temporarily
operate with the assumption that $M=0$.

Taking $N=1$, the above reasoning provides an excellent description,
which may be verified numerically. The concentric semi-circles visible
in the image \ref{cap:Polylog-whole-plane} can be entirely explained
by the behavior of $q^{s}$ as $q\to0$, that is, as $z\to1$. Thus,
the difference between the $N=1$ sheet and the $N=0$ sheet is \begin{equation}
\Delta_{1}=e^{i\pi s/2}\frac{\left(2\pi\right)^{s}}{\Gamma(s)}\left(\frac{\ln z}{2\pi i}\right)^{s-1}\label{eq: N=+1 monodromy}\end{equation}
 Properly speaking, $\Delta_{1}$ is a function of $s$ and $z$,
and so one should write $\Delta_{1}(s;z)$ to signify this. However,
for ease of notation, this marking is dropped, and is taken implicitly
in what follows. The difference $\mbox{Li}_{s}\left(z\right)-\Delta_{1}$
is illustrated in figure \ref{cap:N=3D+1 Upper-sheet}; the concentric
rings are seen to be cancelled out precisely, leaving behind a smoothly
varying function, having the expected smooth structure. Moving across
the cut, for $\Re z>1$, the joint appears to be smooth. For general
positive $N$, it is easy to confirm numerically that the discontinuity
between the $N$'th sheet, and the $N-1$'th sheet, across the $\Re z>1$
cut is \[
\Delta_{N}=e^{i\pi s/2}\frac{\left(2\pi\right)^{s}}{\Gamma(s)}\left(N-1+\frac{\ln z}{2\pi i}\right)^{s-1}\]

The discontinuity across the cut $\Re z>1$, for negative $N$, follows
similarly. Next, taking $N=-1$ in equation \ref{eq:log-discontinuity},
(or $N=0$ in the equation immediately above, which can serve as a
point of confusion) the difference is given by -$1/(-1+q)^{1-s}=e^{-i\pi s}/(1-q)^{1-s}$,
or \[
\Delta_{-1}=e^{-i\pi s/2}\frac{\left(2\pi\right)^{s}}{\Gamma(s)}\left(1-\frac{\ln z}{2\pi i}\right)^{s-1}\]

The figure \ref{cap:N=3D-1 Lower-sheet} illustrates the difference
$\mbox{Li}_{s}\left(z\right)-\Delta_{-1}$. For general negative $N$,
the difference between adjacent sheets is then \[
\Delta_{-N}=e^{-i\pi s/2}\frac{\left(2\pi\right)^{s}}{\Gamma(s)}\left(N-q\right)^{s-1}\]
 Again, it can be numerically verified that that the transition from
sheet to sheet is smooth as one crosses the cut $\Re z>1$.

It is perhaps more clear to use explicit topological language. Let
$m_{1}$ represent the homotopy class of all loops based at some point
$z$ on the complex plane, that wind once around the branch-point
$z=1$ in the positive direction. The action of $m_{1}$ on the polylogarithm
has the effect of carrying the polylog from one sheet to the next.
In the above discussion, it was determined that \[
m_{1}\cdot\mbox{Li}_{s}\left(z\right)=\mbox{Li}_{s}\left(z\right)-\Delta_{1}\]
 The logarithm in $\Delta_{1}$ has a branch point at $z=0$. That
is, after acting with $m_{1}$, one is now on a sheet with a cut extending
from $z=0$ to the right. Let let $m_{0}$ represent the homotopy
class of loops that wind once around the branch point $z=0$ in a
right-handed fashion. Acting on the logarithm, one has \[
m_{0}\cdot\ln z=\ln z+2\pi i\]
 Recalling the definition of $q=\ln z/2\pi i$, one thus has $m_{0}\cdot q=q+1$,
and so \[
m_{0}\cdot\Delta_{N}=\Delta_{N+1}\]
The principal sheet of the polylogarithm has no branch point at $z=0$,
and so one has \[
m_{0}\cdot\mbox{Li}_{s}\left(z\right)=\mbox{Li}_{s}\left(z\right)\]
Winding in the opposite direction, one has \[
m_{1}^{-1}\cdot\mbox{Li}_{s}\left(z\right)=\mbox{Li}_{s}\left(z\right)-\Delta_{-1}\]
 In order for $m_{1}$ to be properly considered as the group-theoretic
inverse of $m_{1}$, one must have $m_{1}\cdot m_{1}^{-1}\cdot\mbox{Li}_{s}\left(z\right)=m_{1}^{-1}\cdot m_{1}\cdot\mbox{Li}_{s}\left(z\right)=\mbox{Li}_{s}\left(z\right)$.
The resolution of this implies $m_{1}\cdot\Delta_{-1}=-\Delta_{1}$,
which in turn implies that $m_{1}\cdot q=q+1$. This seems strange,
as the logarithm has no branch point at $z=1$, and so there is nothing
to wind around. By this argument, one would have expected that $m_{1}$
had no effect on $q$ at all. The point is subtle and is worth establishing
clearly; it is the joining of the polylogarithm to the logarithm cuts
that causes the effect. Consider only the logarithm (not the polylogarithm),
arranged so that the branch extends to the right. Consider starting
just below the real axis, and winding around in a right-handed fashion
around the point $z=1$, and finishing just above the real axis. There
is no obstruction at $z=1$, and so this loop can be shrunk to a very
short line segment. None-the-less, the short line segment crosses
the cut, and its effect is to move to a different sheet. It is for
this reason that one has $m_{1}\cdot q=q+1$. A different set of group
generators will be defined below, that do capture the idea that winding
around $z=1$ should have no effect on $q$.

One can now consider the task of more complex paths from sheet to
sheet. Passing twice around the $z=1$ branch, one has 

\[
m_{1}^{2}\cdot\mbox{Li}_{s}\left(z\right)=m_{1}\cdot\left[\mbox{Li}_{s}\left(z\right)-\Delta_{1}\right]=\mbox{Li}_{s}\left(z\right)-\Delta_{1}-\Delta_{2}\]
 Repeating this gluing $n$ times, so that each time, one pastes the
sheets so that crossing the $z>1$ cut is smooth and differentiable,
gives \begin{eqnarray*}
m_{1}^{n}\cdot\mbox{Li}_{s}\left(z\right) & = & \mbox{Li}_{s}\left(z\right)-\sum_{k=1}^{n}\Delta_{k}\\
 & = & \mbox{Li}_{s}\left(z\right)-e^{i\pi s/2}\frac{\left(2\pi\right)^{s}}{\Gamma(s)}\sum_{k=1}^{n}\frac{1}{\left(k-1+q\right)^{1-s}}\\
 & = & \mbox{Li}_{s}\left(z\right)-e^{i\pi s/2}\frac{\left(2\pi\right)^{s}}{\Gamma(s)}\left[\zeta\left(1-s,\frac{\ln z}{2\pi i}\right)-\zeta\left(1-s,n+\frac{\ln z}{2\pi i}\right)\right]\end{eqnarray*}
 which is recognizable from the initial considerations. 

To capture the idea of there being no obstruction at $z=1$ for the
logarithm, one may define a group element $g_{1}=m_{1}m_{0}^{-1}$,
so that one has $g_{1}\cdot\ln z=\ln z$. Then, define $g_{0}=m_{0}$.
In terms of these generators, one has the relations \begin{eqnarray*}
g_{0}\cdot q & = & q+1\\
g_{1}\cdot q & = & q\\
g_{0}\cdot\mbox{Li}_{s}\left(z\right) & = & \mbox{Li}_{s}\left(z\right)\\
g_{1}\cdot\mbox{Li}_{s}\left(z\right) & = & \mbox{Li}_{s}\left(z\right)-\Delta_{1}\\
g_{0}\cdot\Delta_{N} & = & \Delta_{N+1}\quad\mbox{ for }N>0\\
g_{1}\cdot\Delta_{N} & = & \Delta_{N}\end{eqnarray*}
 To complete the picture, one needs the action of $g_{0}$ on $\Delta_{-1}$,
or equivalently $g_{0}^{-1}$ on $\Delta_{1}$. There is some ambiguity,
and thus, room for confusion, as the term $(-1)^{s}$ arises, which
may be taken as $e^{i\pi s}$ or as $e^{-i\pi s}$, with the last
two inequivalent for non-integer $s$. The resolution lies in considering
the joining of sheets across the cut that runs between $0<\Re z<1$.
Visually, the cut can be clearly seen in figure \ref{cap:N=3D-1 Lower-sheet}.
Consider now the task of glueing sheets across this cut. For a point
$z$ just above the line connecting $z=0$ and $z=1$, one has $q=\epsilon+iv$
for some small, positive $\epsilon$ and some general, positive $v$.
Just below this line, one has $q=1-\epsilon+iv$. That is, $q$ differs
by 1, of course. This cut cannot be glued to the polylog, of course;
it must be glued to another sheet of the logarithm. The correct gluing,
for $z=x$ real, $0<x<1$, is given by \[
\lim_{\epsilon\to0}\left[\Delta_{1}(x+i\epsilon)+e^{i2\pi s}\Delta_{-1}(x-i\epsilon)\right]=0\]
 That is, \[
g_{0}\cdot\Delta_{-1}=-e^{-i2\pi s}\Delta_{1}\]
 This somewhat strange form is nothing more than an accounting trick;
it is the result of providing a definition of $\Delta_{-N}$ that
had a {}``natural'' normalization, avoiding a minus sign. The price
of that definition is this glitch. In all other respects, the homotopy
proceeds as expected, so that, for example, \[
g_{0}\cdot\Delta_{-N}=\Delta_{-N+1}\]
 for $N>1$. 

The free combinations of powers of the two operators $g_{0}$ and
$g_{1}$ (or $m_{0}$ and $m_{1}$) generate a group, the monodromy
group of the polylogarithm. For $s=m$ a positive integer, the monodromy
group has a finite-dimensional representation of dimension $m+1$.
Well-known is the case for $s=2$, the dilogarithm, where the representation
forms the discrete Heisenberg group. This may be seen as follows.
For $s=2,$ one has \[
\Delta_{n}=2\pi i\left(\ln z+(n-1)2\pi i\right)\]
 Repeated applications of $g_{0}$ and $g_{1}$ result in a linear
combinations of $\mbox{Li}_{2}\left(z\right)$, $\ln z$ and 1; no
terms of any other type appear. One may take these three elements
as the basis of a three-dimensional vector space. Writing the basis
as column vectors, the representation is \begin{eqnarray*}
4\pi^{2}\mapsto e_{1} & = & \left[\begin{array}{c}
1\\
0\\
0\end{array}\right]\\
-2\pi i\ln z\mapsto e_{2} & = & \left[\begin{array}{c}
0\\
1\\
0\end{array}\right]\\
\mbox{Li}_{2}\left(z\right)\mapsto e_{3} & = & \left[\begin{array}{c}
0\\
0\\
1\end{array}\right]\end{eqnarray*}
 so that the action may be represented as \[
g_{1}=\left[\begin{array}{ccc}
1 & 0 & 0\\
0 & 1 & 1\\
0 & 0 & 1\end{array}\right]\quad\mbox{ and }\quad g_{0}=\left[\begin{array}{ccc}
1 & 1 & 0\\
0 & 1 & 0\\
0 & 0 & 1\end{array}\right]\]
 These two matrices may be seen to be the generators of the discrete
Heisenberg group $\mathcal{H}_{3}(\mathbb{Z})$. The general group
element of the discrete Heisenberg group is \[
\left[\begin{array}{ccc}
1 & a & b\\
0 & 1 & c\\
0 & 0 & 1\end{array}\right]\]
 for integers $a,b,c$. By convention, one writes $x$ for $g_{0}$
and $y$ for $g_{1}$, and defines a new element $z$ (having no relation
at all to the argument $z$ of the dilogarithm) as the commutator
$z=xyx^{-1}y^{-1}$, so that \[
z=\left[\begin{array}{ccc}
1 & 0 & 1\\
0 & 1 & 0\\
0 & 0 & 1\end{array}\right]\]
 The element $z$ is in the center, so that it commutes with $x$
and $y$, so $zx=xz$ and $zy=yz$. These relations may be taken as
giving the group presentation, so that \[
\mathcal{H}_{3}(\mathbb{Z})=\left\langle (x,y)\,\vert\; z=xyx^{-1}y^{-1},\; zx=xz,\; zy=yz\right\rangle \]
 Curiously, one should note that the discrete Heisenberg group has
an alternate presentation that is reminiscent of the braid group.
This is perhaps not a complete surprise, as a monodromy can be generally
understood to be a subgroup of the braid group or mapping class group
of the punctured disk. The mapping class group permutes the punctures
of the disk in a continuous fashion, so that the result is an intertwining
of strands. For a disk with two punctures, the braid group is $B_{3}$
consists of three strands. In this case, two of the strands can be
visualized as two parallel lines, each representing one of the two
branch points of the polylogarithm. The third strand is the monodromy
path, which weaves about the two straight lines, but must always return
to the middle. Now, the braid group $B_{3}$ has two group generators,
denoted $\sigma_{1}$ and $\sigma_{2}$ by convention, which denote
a single twist of the left or right pair of strands of the braid.
To weave the middle strands of the monodromy about the two branch
points, one identifies $g_{0}\to\sigma_{1}^{2}$ and $g_{1}\to\sigma_{2}^{2}$.
At this point, ignoring the polylog, and limiting ones attention to
the universal covering space for the homotopy group over two branch
points, one very simply has that $g_{0}$ and $g_{1}$ are the generators
of the free group in two letters. Now for the curious resemblance.
The braid group $B_{3}$ has the group presentation $\sigma_{1}\sigma_{2}\sigma_{1}=\sigma_{2}\sigma_{1}\sigma_{2}$.
The discrete Heisenberg group has the remarkably similar identity
$xyyx=yxxy$.

For $s=3$, one may write \[
g_{0}=\left[\begin{array}{cccc}
1 & 1 & 1 & 0\\
0 & 1 & 2 & 0\\
0 & 0 & 1 & 0\\
0 & 0 & 0 & 1\end{array}\right]\quad\mbox{ and }\quad g_{1}=\left[\begin{array}{cccc}
1 & 0 & 0 & 0\\
0 & 1 & 0 & 0\\
0 & 0 & 1 & 1\\
0 & 0 & 0 & 1\end{array}\right]\]
 for the basis \begin{eqnarray*}
4\pi^{3}i\mapsto e_{1} & = & \left[\begin{array}{c}
1\\
0\\
0\\
0\end{array}\right]\\
2\pi^{2}\ln z\mapsto e_{2} & = & \left[\begin{array}{c}
0\\
1\\
0\\
0\end{array}\right]\\
-\pi i\left(\ln z\right)^{2}\mapsto e_{3} & = & \left[\begin{array}{c}
0\\
0\\
1\\
0\end{array}\right]\\
\mbox{Li}_{3}\left(z\right)\mapsto e_{4} & = & \left[\begin{array}{c}
0\\
0\\
0\\
1\end{array}\right]\end{eqnarray*}
The generated group is not the four-dimensional Heisenberg group,
as there simply aren't enough generators for that. Its not a representation
of the three-dimensional Heisenberg group, as the commutator $w=g_{0}g_{1}g_{0}^{-1}g_{1}^{-1}$
is not in the center, since $wg_{0}\ne g_{0}w$. However, one does
have $wg_{1}=g_{1}w$, and from this one may deduce that the monodromy
group has the presentation \begin{equation}
M=\langle g_{0},g_{1}\,\vert\, g_{1}g_{0}g_{1}^{-1}g_{0}^{-1}=g_{0}g_{1}^{-1}g_{0}^{-1}g_{1}\rangle\label{eq:polylog monodromy}\end{equation}
 In any case, the group is solvable and thus thin. 

The Heisenberg group may be regained as the quotient group $\mathcal{H}_{3}(\mathbb{Z})=M/\left\langle w\right\rangle $
where $\left\langle w\right\rangle $ is the conjugate closure of
$w$: \[
\left\langle w\right\rangle =\left\{ gwg^{-1}\,\vert\, g\in M\right\} \]
 The conjugate closure is, of course, a normal subgroup of $M$. Since
$g_{1}$ already commuted with $w$, what the quotient group construction
does is to impose $\left\langle w\right\rangle g_{0}=g_{0}\left\langle w\right\rangle $,
which shows up on the cosets as the other, {}``missing'' relationship
needed to get the Heisenberg group. It is thus that the presentation
of the Heisenberg group is regained. The general form of an element
$h\in\left\langle w\right\rangle $ may be taken to be \[
h=g_{0}^{-n}wg_{0}^{n}=\left[\begin{array}{cccc}
1 & 0 & 0 & 2n+1\\
0 & 1 & 0 & -2\\
0 & 0 & 1 & 0\\
0 & 0 & 0 & 1\end{array}\right]\]
 for integer $n$. Thus, it is clear that $\left\langle w\right\rangle $
is abelian, and that $\left\langle w\right\rangle \cong\mathbb{Z}$.

For general integer $s=m$, the monodromy group is always unipotent,
and thus a nilpotent group. The generators take the form \[
g_{0}=\left[\begin{array}{cc}
C & \begin{array}{c}
0\\
\vdots\\
0\end{array}\\
0\;\cdots\;0 & 1\end{array}\right]\quad\mbox{ and }\quad g_{1}=\left[\begin{array}{ccccc}
1 & 0 & \cdots &  & 0\\
0 & 1 &  &  & \vdots\\
\vdots &  & \ddots & 0 & 0\\
 &  & 0 & 1 & 1\\
0 & \cdots & 0 & 0 & 1\end{array}\right]\]
 where $C$ are the binomial coefficients written in an upper-triangular
fashion (thus, for example, for $m=4$, one adds a column of $(1,3,3,1)$
and, for $m=5$, another column $(1,4,6,4,1)$, \emph{etc}.). This
form is determined by taking as the basis vectors the monomials appearing
in the expansion of $\left(\ln z+2\pi i\right)^{m}$ (whence the binomial
coefficients), and normalizing with a factor of $2\pi i/\Gamma(m)$.
The generated monodromy group is isomorphic to the above $s=3$ case,
as it has the same presentation. The element $w=g_{0}g_{1}g_{0}^{-1}g_{1}^{-1}$
commutes with $g_{1}$, just as before; but $g_{0}w\ne wg_{0}$. Thus,
the conjugate closure $\left\langle w\right\rangle $ is again a one-dimensional
abelian group, that is, $\left\langle w\right\rangle \cong\mathbb{Z}$,
and the resulting quotient group is again the Heisenberg group.

Note that the existence of polylogarithm ladders seems to be equivalent
to the statement that products of monodromy groups can be factored
in interesting and unusual ways. XXX Expand on this remark.

For $s$ not an integer, the action of $g_{0}$ does not close, and
the general vector-space representation is infinite dimensional. The
action of the $g_{0}$ monodromy can be understood to be the shift
operator on an infinite-dimensional vector space, whose basis vectors
may be taken to be $\Delta_{N}$; the operator $g_{1}$ acts to inject
into this shift space. Defining $w$ as before, it is not hard to
determine that $w\cdot\Delta_{N}=\Delta_{N}$ and that $w\cdot\mbox{Li}_{s}\left(z\right)=\mbox{Li}_{s}\left(z\right)+\Delta_{1}-\Delta_{2}$.
Remarkably, the same commutation relation holds as before, in that
$g_{1}w\cdot\mbox{Li}_{s}\left(z\right)=wg_{1}\cdot\mbox{Li}_{s}\left(z\right)=\mbox{Li}_{s}\left(z\right)-\Delta_{2}$.
This, together with the more trivial result $g_{1}w\cdot\Delta_{N}=wg_{1}\cdot\Delta_{N}=\Delta_{N}$
shows that the infinite-dimensional vector space representation of
the monodromy has the same presentation as the finite-dimensional
case. Thus, one concludes that, for all complex values of $s\ne2$,
the monodromy of the polylogarithm is given by equation \ref{eq:polylog monodromy}. 

It is not clear whether there may also be smaller, finite-dimensional
representations for particular values of $z$. XXX Expand on this.

The analog of the Dirichlet $L$-functions are functions of the form
\[
p^{-s}\sum_{m=1}^{p-1}\chi(m)\,\mbox{Li}_{s}\left(ze^{2\pi im/p}\right)\]
 where $\chi(m)$ is a primitive character modulo $p$. These functions
have the strange property of having branch points at $e^{2\pi im/p}$
whenever $\chi(m)$ is not zero. The resulting monodromy groups have
a more complex structure as well.

\section{The {}``Periodic'' Zeta Function }

\begin{figure}

\caption{\label{cap:Periodic-zeta-function}Periodic zeta function}

\includegraphics[%
  width=1\textwidth]{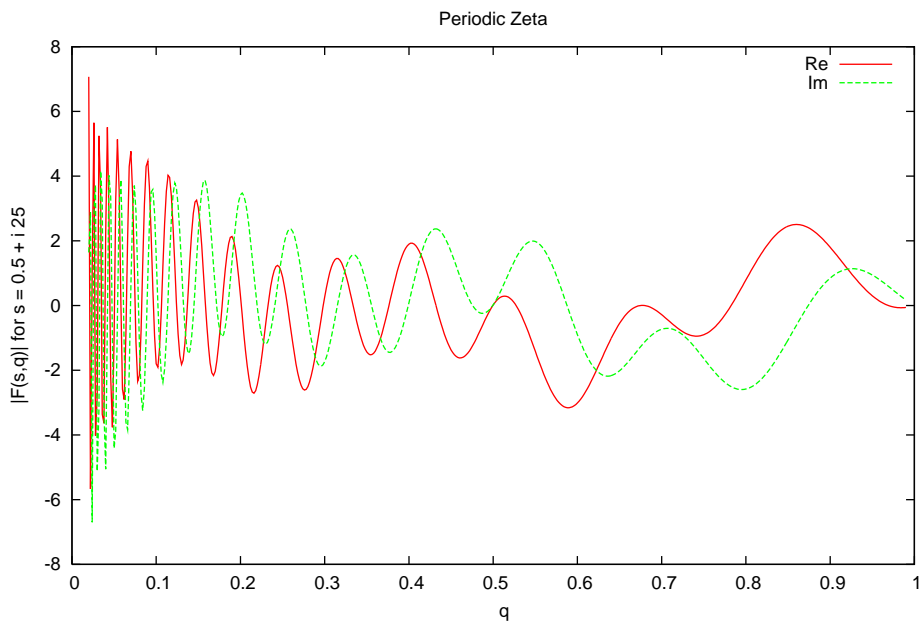}

This graph shows the real and imaginary parts of the periodic zeta
function \[
F(q,s)=\mbox{Li}_{s}\left(e^{2\pi iq}\right)\]
 for $s=\frac{1}{2}+i25$. This value of $s$ is close to a zero of
the Riemann zeta function, at $s=\frac{1}{2}+i25.01085758014\ldots$
Thus, both the real and imaginary parts approach zero at $q=1$, as
well as at $q=\frac{1}{2}$. The increasing oscillations as $q\to0$
are due to the contribution of the very first term of the Hurwitz
zeta: that is, these oscillations are nothing other than that of $q^{-s}=e^{i25\log q}/\sqrt{q}$.
Subtracting these away is the same as analytically continuing to the
region $q>1$, and matches the coarser set of oscillations, which
are given by $(1+q)^{-s}=e^{i25\log(1+q)}/\sqrt{1+q}$

Noteworthy is that the presumed {}``periodicity'' in $q$ is very
misleading: the image suggests an essential singularity at $q=0$,
and continuing, logarithmically slower oscillatory behavior for the
analytic continuation to the region where $q>1$.
\end{figure}

The ability to compute the polylogarithm allows one to visualize some
of the interesting sums that occur in number theory. One such example
is the so-called {}``periodic zeta function'', defined by Apostol\cite[Thm 12.6]{Apo76}
as \begin{equation}
F(q;s)=\sum_{n=1}^{\infty}\frac{e^{i2\pi nq}}{n^{s}}\label{eq:periodic zeta}\end{equation}
 Clearly, one has $F(q,s)=\mbox{Li}_{s}\left(e^{2\pi iq}\right)$.
The periodic zeta function occurs when one attempts to state a reflection
formula for the Hurwitz zeta function, as \begin{equation}
\zeta(1-s,q)=\frac{\Gamma(s)}{(2\pi)^{s}}\left[e^{-i\pi s/2}F(q;s)+e^{i\pi s/2}F(1-q;s)\right]\label{eq:Hurwitz reflection}\end{equation}
 What makes the periodic zeta function interesting is that \emph{it
is not actually periodic}. That one might think it is seems obvious
from the definition \ref{eq:periodic zeta}: the direct substitution
of $q\to q+1$ gives the false suggestion of periodicity in $q$.
This is false because, in fact, $\mbox{Li}_{s}\left(z\right)$ has
a branch point at $z=1$. The {}``periodic'' zeta function is multi-sheeted,
and, attempting to trace a continuous path from $q\to q+1$, while
keeping $q$ real carries one through the branch-point, and from one
sheet to the next. This is illustrated graphically, in figure \ref{cap:Periodic-zeta-function}.

The figure \ref{cap:Periodic-zeta-function} shows an oscillatory
behavior that clearly diverges as $q\to0$. This oscillation is can
be simply explained. Substituting $s=\frac{1}{2}+i\tau$ into equation
\ref{eq:Hurwitz reflection} gives \[
\zeta\left(\frac{1}{2}-i\tau,\, q\right)=\frac{e^{-i\pi/4}e^{-i\tau\log2\pi}e^{i\phi}}{\sqrt{2\cosh\tau\pi}}\left[e^{\pi\tau/2}F(q;s)+ie^{-\pi\tau/2}F(1-q;s)\right]\]
 after the substitution of \[
\Gamma\left(\frac{1}{2}+i\tau\right)=\frac{\sqrt{\pi}e^{i\phi}}{\sqrt{\cosh\tau\pi}}\]
 where $\phi$ is a purely real phase that can be expressed as a somewhat
complicated sum\cite[eqn 6.1.27]{A&S}. For reasonably large, positive
$\tau$, such as $\tau=25$ in the picture, one may ignore the second
term, so that the whole expression simplifies to \[
F\left(q;\frac{1}{2}+i\tau\right)=\exp\, i\left(\frac{\pi}{4}-\phi+\tau\log2\pi\right)\;\zeta\left(\frac{1}{2}-i\tau,\, q\right)+\mathcal{O}\left(C^{-\tau}\right)\]
 for some real constant $C>0$. Then, as $q\to0$, the leading contribution
to the Hurwitz zeta comes from the $n=0$ term in equation \ref{eq:Hurwitz zeta}:
that is, $q^{-s}=e^{i\tau\log q}/\sqrt{q}$, so that \begin{equation}
F\left(q;\frac{1}{2}+i\tau\right)\approx\frac{e^{i\tau\log q}e^{i\psi}}{\sqrt{q}}\;\mbox{ as }q\to0\label{eq:F-approx}\end{equation}
 for some fixed, real phase $\psi$ that is independent of $q$. As
is eminently clear from both the picture \ref{cap:Periodic-zeta-function},
and from the approximation \ref{eq:F-approx}, the limit of $q\to0$
of $F(q,s)$ cannot be taken: this is the branch-point of $\mbox{Li}_{s}\left(z\right)$
at $z=1$.

Curiously, the limit $q\to1$ does exist, and one has $F(1,s)=\zeta(s)$
the Riemann zeta. The approximation \ref{eq:F-approx} hints at how
to move through the branch point, from one sheet to the next: \ref{eq:F-approx}
is the discontinuity between sheets. That is, one has, for large,
positive $\tau$\[
F\left(q+1;\,\frac{1}{2}+i\tau\right)\approx F\left(q,\,\frac{1}{2}+i\tau\right)-\frac{e^{i\tau\log q}e^{i\psi}}{\sqrt{q}}\]
 as the formula that approximates the movement from one sheet to the
next. In essence, this shows that the {}``periodic'' zeta function
is not at all periodic: $F(q+1;s)\ne F(q;s)$ whenever $\Re s\le1$. 

The complete relationship between the Hurwitz zeta and the periodic
zeta is rather subtle. For example, if instead one considers large
\emph{negative} $\tau$, and graphs the periodic zeta, one obtains
what is essentially the left-right reversed image of \ref{cap:Periodic-zeta-function},
with oscillations approaching a singularity at $q\to1$. One may easily
\emph{numerically} verify that these oscillations are precisely of
the form $(1-q)^{-s}$. From this, one may deduce that the correct
form of the reflection formula is \begin{equation}
F(q,1-s)=\frac{\Gamma(s)}{(2\pi)^{s}}\left[e^{i\pi s/2}\zeta(s,q)+e^{-i\pi s/2}\zeta(s,1-q)\right]\label{eq:periodic reflection}\end{equation}
 which captures the oscillatory behavior at $q\to0$ for large positive
$\tau$ and the oscillations at $q\to1$ for large negative $\tau$. 

The constraint of keeping $q$ real causes integer values of $q$
to correspond precisely to the branch point of the polylogarithm.
This makes reasoning about the continuity of the periodic zeta at
integral values of $q$ rather turbid. Since the Riemann zeta function
lies precisely at the branch point, this seems to also make it difficult
to gain new insight into the Riemann zeta in this way.

\section{Conclusion}

Both the Borwein and the Euler-Maclaurin algorithms appears to offer
a stable and fast way of computing the Hurwitz zeta function for general
complex $s$ and real $q$, and the polylogarithm for general complex
$s$ and $z$. The Euler-Maclaurin algorithm offers superior performance
for the Hurwitz zeta. An actual implementation using a variable-precision
library shows that all of the algorithms are quite tractable. 

An unexplored area is the optimal implementation of the algorithms
using standard IEEE double-precision floating-point mathematics. Preliminary
work shows that some of the intermediate terms are just large enough
so that rounding error may be of concern; the region of high-precision
convergence has not been characterized. The propagation of rounding
errors through the algorithm has not been characterized; it may be
possible to re-arrange terms to minimize the propagation of errors.
Such a characterization is necessary before reliable double-precision
code can be created. By contrast, such considerations can be swept
under the rug when employing arbitrary-precision arithmetic.

For the Borwein algorithm, evaluating the bounds \ref{eq:error-bound-positive}
and \ref{eq:error-bound-negative}, so as to obtain an order estimate,
is non-trivial, and can account for a significant amount of compute
time. A speedier heuristic is desirable. By contrast, the Euler-Maclaurin
algorithm seems to have a very simple heuristic that is quite adequate
for general-purpose applications. In either case, it would be desirable
to have a more refined analysis and presentation of the heuristics.

It may be possible to find even more rapidly converging algorithms,
by generalizing some of the results from Cohen \emph{et al}\cite{Coh00},
or by appealing to the general theory of Padé approximants. Both approaches
seem promising, although that of Padé approximants may yield a more
general theory.

The polylogarithm and Hurwitz zeta functions are both special cases
of the Lerch transcendental \[
\sum_{n=0}^{\infty}\frac{z^{n}}{(n+q)^{s}}\]
 It should be possible to extend the techniques in this paper to provide
a rapid, globally convergent algorithm for the Lerch transcendental,
thus enabling a deeper numerical exploration of its peculiarities. 

This paper concludes with the application of the algorithm to render
some pretty images.

\bibliographystyle{amsplain}
\bibliography{/home/linas/src/fractal/paper/fractal}

\appendix

\section*{Appendix: Some Pretty Pictures}

Below follow some interesting pictures of the polylogarithm and the
Hurwitz zeta function for values of $s$ on the critical line $s=1/2+i\tau$.
They are shown here, as they are not commonly known. Each figure represents
hours of computation on modern computers. 

\begin{figure}

\caption{The Periodic Zeta}

\includegraphics[%
  width=1\textwidth]{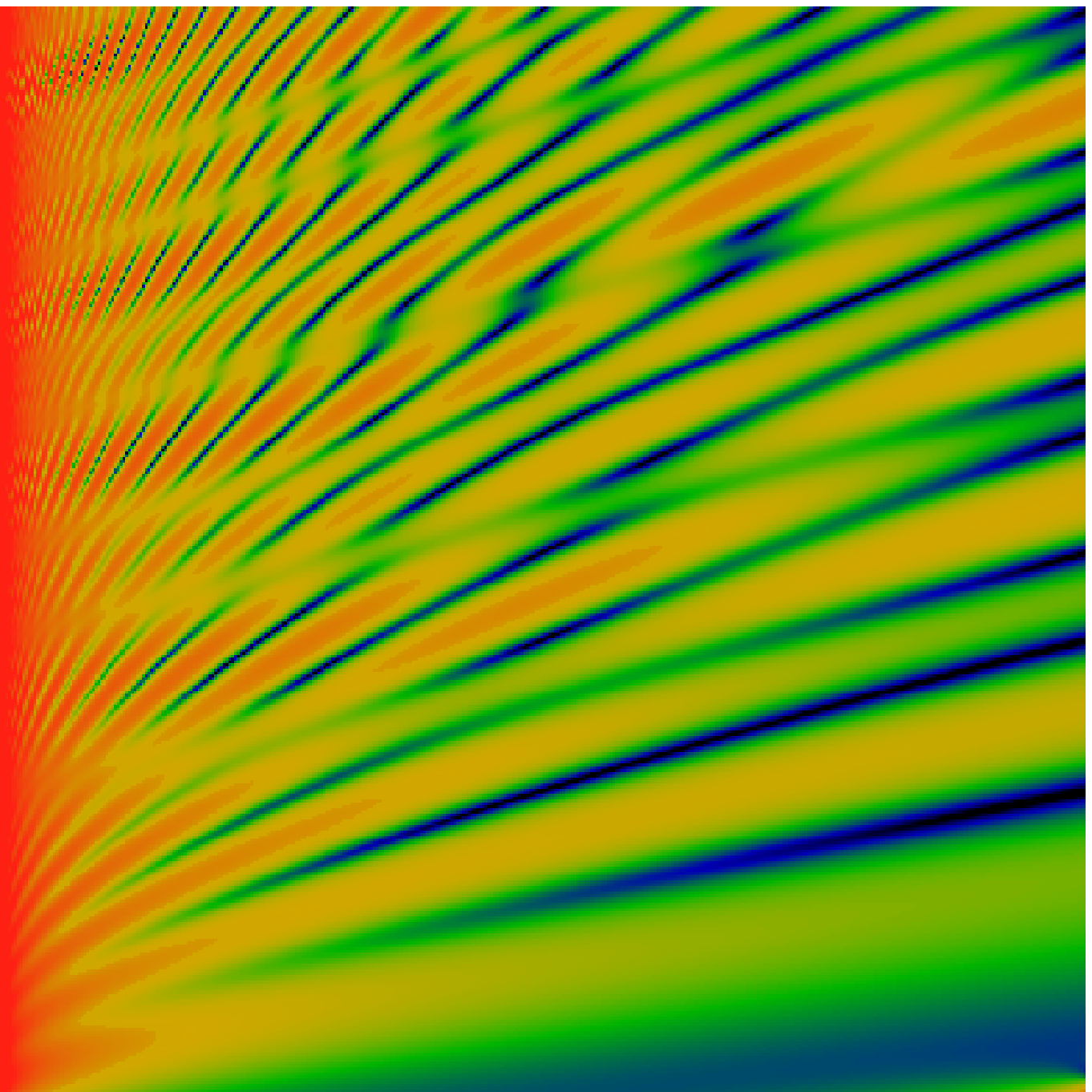}

This image illustrates the so-called {}``periodic zeta function''
\[
F(q;s)=\mbox{Li}_{s}\left(e^{2\pi iq}\right)=\sum_{n=1}^{\infty}\frac{\exp(2i\pi nq)}{n^{s}}\]
 Graphed is the magnitude of $F(q;s)$, with a color scale such that
black represents zero, greenish-blue being values on the order of
one, yellow on the order of two, with increasingly large values shown
as orange-red. Along the horizontal axis runs $q$, taken to be real,
from $q=0$ on the left to $q=1$ on the right. Then, writing $s=\frac{1}{2}+i\tau$,
the value of $\tau$ is varied along the vertical axis, running from
$\tau=0$ at the bottom of the image, to $\tau=50$ at the top of
the image. Zeroes of the Riemann zeta function $\zeta(s)=F(1,s)$
are located where the blue lines intersect the right edge of the image.
From the bottom, the first three zeroes are at $s=0.5+i14.13,\,0.5+i20.02,\,0.5+i25.01$.
Due to the relation to the Dirichlet eta function, the zeros also
materialize at the same values of $s$, but on the $q=1/2$ line.

Remarkably, the blue streaks seem to be roughly parabolic, but are
interrupted by nearly straight {}``erasures''. The pattern is reminiscent
of graphs of the strong-field Stark effect (need ref). In the Stark
effect, eigenvalues are given by the characteristic values of the
Mathieu functions. These cross over one another in a curious fashion;
see for example, figure 20.1 in Abramowitz \& Stegun. 

The precise form of the parabolic blue streaks are due to a term of
the form $q^{-s}$, as given by equation \ref{eq:F-approx}. These
may be subtracted, as illustrated in the next image.

\end{figure}

\begin{figure}

\caption{Hurwitz zeta, with leading term subtracted}

\includegraphics[%
  width=1\textwidth]{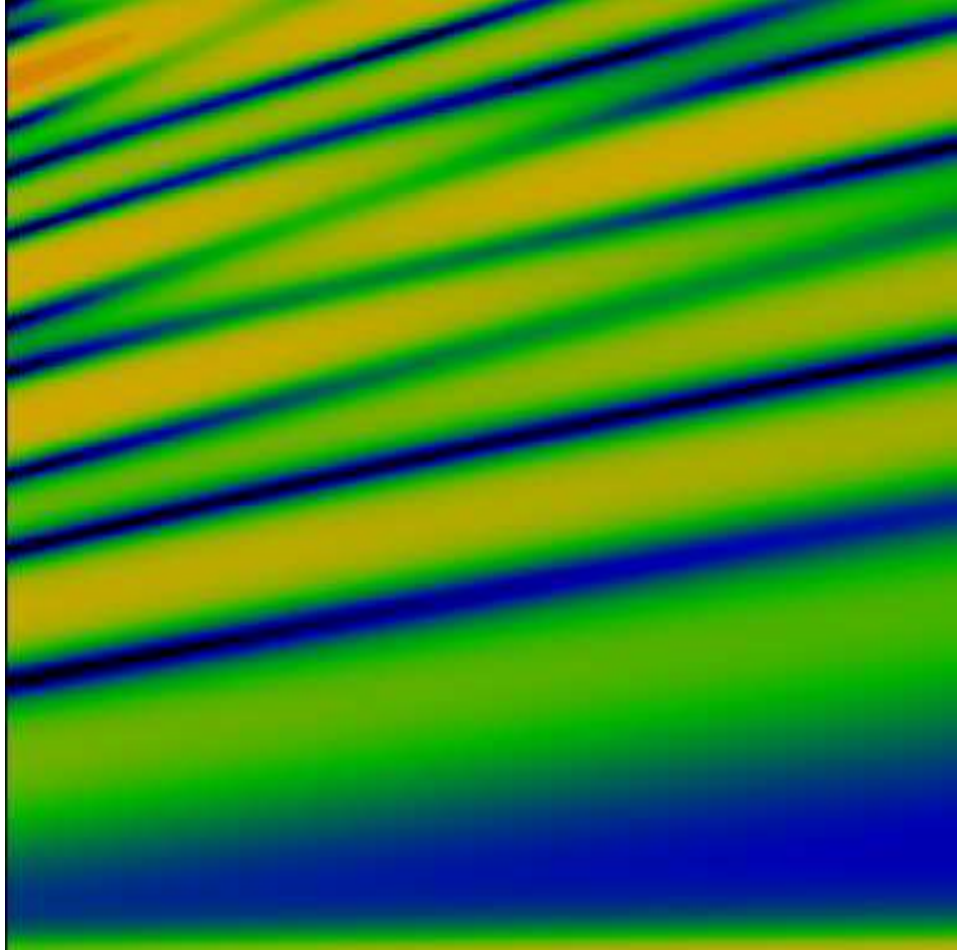}

This image shows the magnitude of \[
\zeta(s,q+1)=\zeta(s,q)-q^{-s}\]
 for the range of $0\le q\le1$along the horizontal axis, and $0\le\tau\le50$
along the vertical axis, just as in the previous image; the same coloration
scheme is used as well. This essentially corresponds to the previous
picture, with the leading divergence subtracted; for large positive
$\tau$, the magnitude of the periodic zeta and the Hurwitz zeta differ
by an exponentially small amount. A careful comparison of this picture
to the previous now shows that the {}``erasures'' or {}``streaks''
in the previous image correspond exactly to the dark parts of this
image. Again, a criss-cross lattice pattern can be seen. The dominant
streaks in this picture are presumably due to the $(1+q)^{-s}$ term. 
\end{figure}

\begin{figure}

\caption{Extended view of Hurwitz zeta}

\includegraphics[%
  width=1\textwidth]{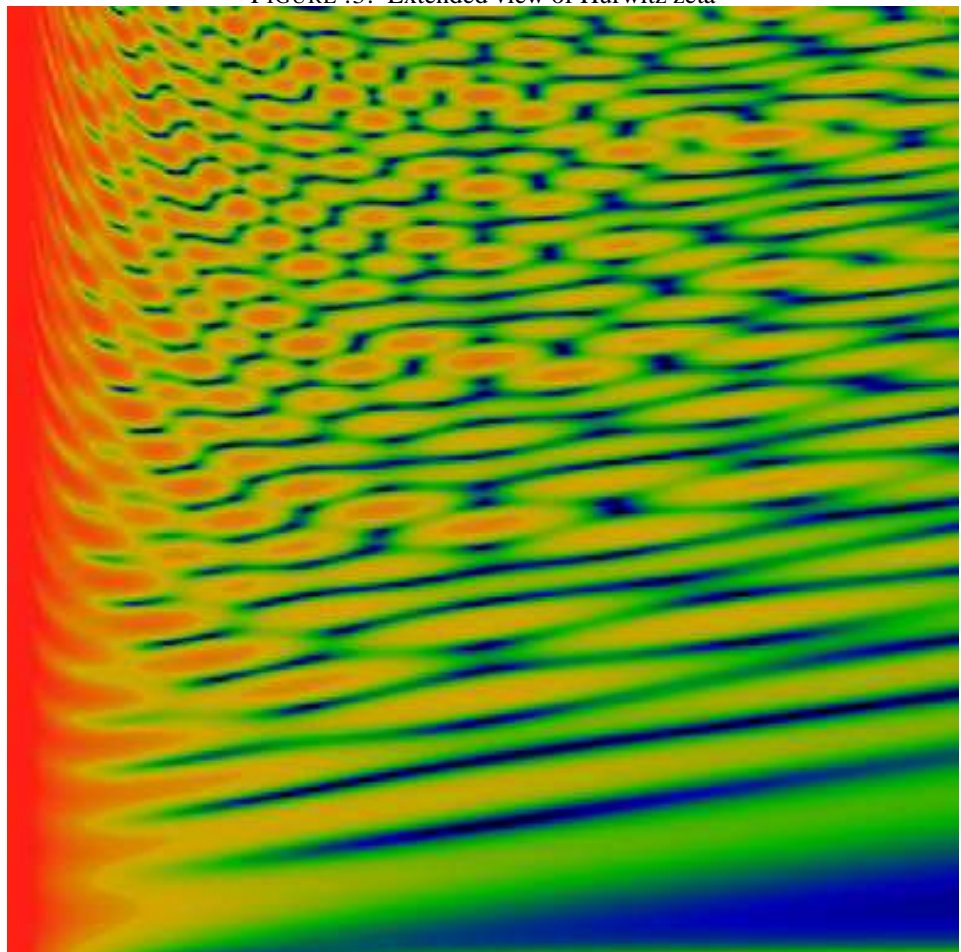}

This figure shows a view of the magnitude of the Hurwitz zeta $\left|\zeta(s,q)\right|$
over an extended range, together with a rescaling of coordinates in
an attempt to straighten out the striations. The value of $q$ ranges
from 0 to $2\sqrt{2}\approx2.8$ along the horizontal axis, although
it is scaled as $q^{3/2}$ along this axis. This rescaling essentially
turns the parabolic striations into straight, radial lines emanating
from the origin at the lower left. This image also re-scales the $\tau$
coordinate, in an attempt to turn the radial striations into horizontals.

The value of $q=1$ is a vertical line running exactly down the middle
of the image; the non-trivial Riemann zeroes are visibly stacked up
on this line. The value of $\tau$ runs from 0 to 100 on this vertical
line. On other verticals, $\tau$ has been rescaled by $q^{1/2}$,
so that $\tau$ runs from 0 to 141 on the right hand side of the image,
while being essentially held at zero along the left edge of the image.

The increasing complexity of the chain-link pattern is apparent as
one moves up the imaginary axis. Equally apparent is that the maximum
range of excursion of the magnitude changes only slowly over the range
of the image: as the pattern increases in complexity, it does not
wash out, but has a distinct binary on/off character. 
\end{figure}

\begin{figure}

\caption{\label{cap:Polylogarithm-phase-plot}Polylogarithm phase plot}

\includegraphics[%
  width=1\textwidth]{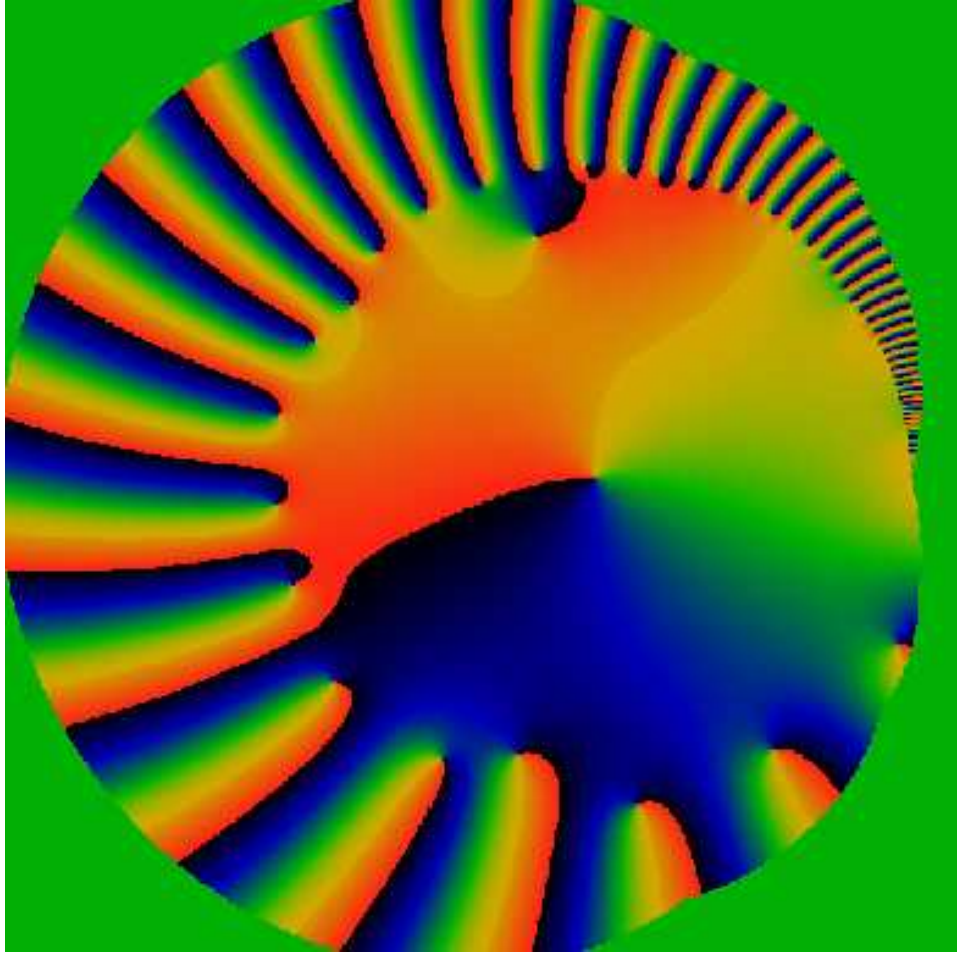}

This image shows a phase plot for the polylogarithm, on the complex
$z$-plane, for $s=0.5+i80$. The color scheme is such that black
represents the phase $\arg\mbox{Li}_{s}\left(z\right)=-\pi$, red
represents the value of $+\pi$, and green a phase of zero. Points
where the phase wraps around the full color range are zeros. Thus,
each of the {}``fingers'' in this image terminates on a zero. Most
of these zeros are close to, but usually not quite on the circle $\left|z\right|=1$.
The lone zero near the center of the image corresponds to the zero
at $z=0$. The slight cusp at the right side of the image corresponds
to $z=1$. The outer boundary of the image, adjoining to the solid
green area, represents the limit of convergence for the combination
of the Borwein algorithm and the duplication formula; one can enlarge
the region of convergence slightly, but not much, but using polynomial
approximations of higher order. The boundary is easily crossed by
applying the inversion formula, as the later images show. The image
is offset somewhat, rather than being centered on $z=0$, because
most of the convergent region is to the left of the imaginary axis.

Note the other lone zero floating along inside of the $\left|z\right|<1$
disk. If this picture were animated as $\tau$ increased from 0 in
the positive direction, then the zero fingers can be see be seen as
marching across in a counter-clockwise fashion, starting at $z=1$and
proceeding around. Most zeros pass to another sheet upon crossing
$z=1$, after making a full revolution; others spiral around a second
time on this sheet, such as the lone zero up top. An animated movie
of this image, showing the spiraling, can be seen at http://www.linas.org/art-gallery/polylog/polylog.html

\end{figure}

\begin{figure}

\caption{\label{cap:Polylog-whole-plane}Polylog whole plane}

\includegraphics[%
  width=1\textwidth]{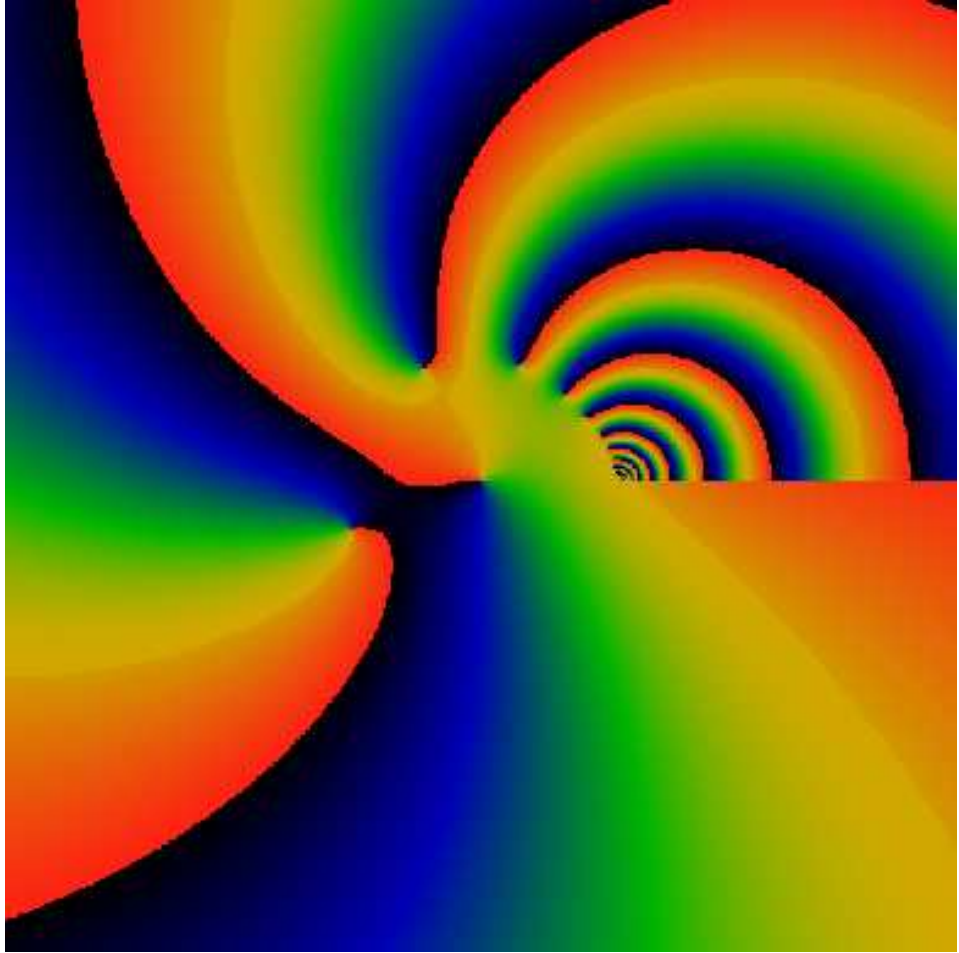} 

This image illustrates the results of combining the Borwein algorithm
together with the inversion and duplication formulas. The image shows
the complex $z$-plane, is centered on $z=0,$ and runs from -3.5
to +3.5, left to right, top to bottom. The value of $s$ is fixed
at $s=0.5+15i$ for the whole picture. The color scheme is the same
as for the previous picture. 

Besides the zero at $z=0$ at the center of the picture, another zero
is visible to the left, and slightly down. This zero is an {}``ex-Riemann
zero'', in that, if instead one had created the image for $s=0.5+i14.1347\ldots$,
then this point would have been located exactly at $z=-1$. As it
happens, this point has now rotated to slightly below $z=-1$. The
zeroes above and to the right of the origin are Riemann zeroes to
be: as $\tau$ increases, each will in turn rotate counterclockwise
to assume the position at $z=-1$. 

The concentric shells are centered on the point $z=1$. The branch
cut can be seen extending to the right from from $z=1$, on the positive
real axis. 
\end{figure}

\begin{figure}

\caption{Polylogarithm off the critical axis}

\includegraphics[%
  width=1\textwidth]{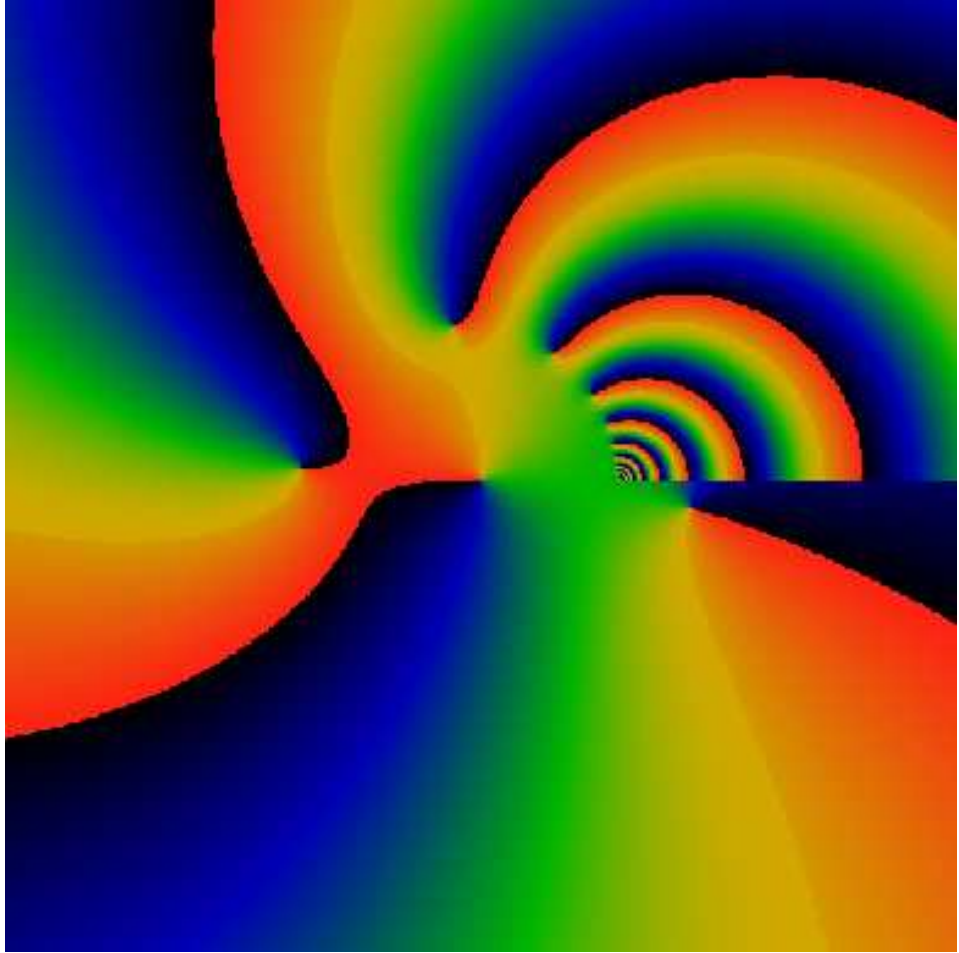}

This figure shows the polylogarithm in the complex $z$-plane, at
fixed $s=\sigma+i\tau=1.2+i14$. It is superficially similar to the
previous image, with two notable differences. Visible just under the
cut, on the right, is a zero. If $\sigma$ were smaller, this zero
would move to the left; at $\sigma=0.5$, it would be very near to
the branch point at $z=1$. As $\tau$ gets larger, this point moves
up, crossing the branch cut and moving onto the next sheet. Of course,
if instead one had $\sigma+i\tau=0.5+i14.13\ldots$, this zero would
be located precisely at $z=1$. Thus, this figure illustrates what
a Riemann zero looks like when it is not on the critical strip.

As $\tau$ increases, the remaining zeroes circle around in a widening
spiral: thus, the zero at the far left of the image is not near $z=-1$,
but is to the left of there (and thus is not going to be a zero of
the Dirichlet eta function). If instead one had $\sigma$ with a value
of less than 0.5, the zeroes would spiral around in an ever-tightening
spiral, and would never hit the branch cut or cross over. The Riemann
hypothesis is essentially that these zeroes are impacting exactly
at the branch point; a violation of RH would be a pair of zeroes that
failed to hit the branch point, instead passing to the left and right
of the branch point. This figure suggests that RH holds: the zeroes
are clearly lined up in single file; it is hard to imagine how this
single file might break ranks into a pair of zeros that simultaneously
passed near the branch point.
\end{figure}

\begin{figure}

\caption{\label{cap:N=3D+1 Upper-sheet}Upper sheet of the polylogarithm}

\includegraphics[%
  width=1\textwidth]{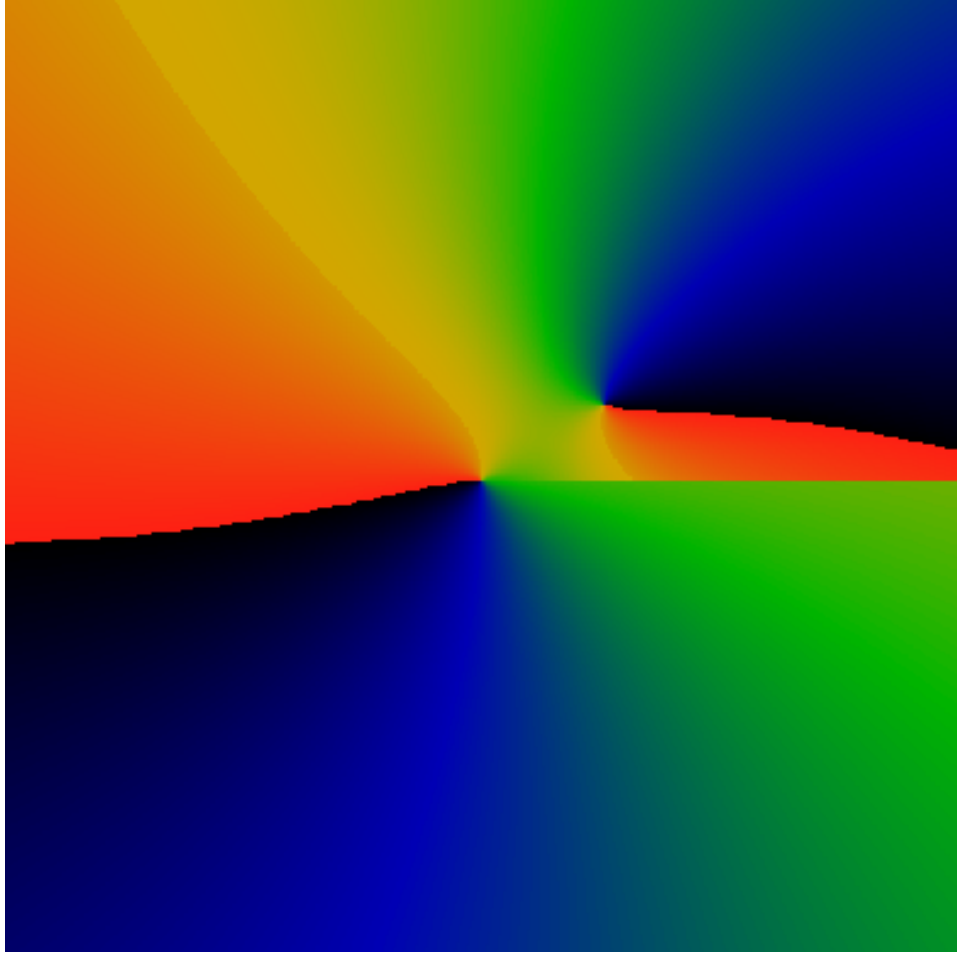}

This figure shows the $g_{1}\mbox{Li}_{s}(z)=\mbox{Li}_{s}(z)-\Delta_{1}$
sheet of the polylogarithm on the complex $z$-plane, obtained by
circling the branch-point at $z=1$ in the right-handed direction.
Two zeroes are visible: one at $z=0$, and the other an {}``ex-Riemann
zero'', located not far from $z=1$, a bit above the real axis. For
this image, $s$ is held constant, at $s=0.5+i15$; the Riemann zeta
function has its first non-trivial zero at $s=0.5+i14.13\ldots$.
At that value, the above zero would have been located precisely at
$z=1$, on the branch point of the polylogarithm. As $\tau$ increases,
the Riemann zeroes pass precisely through this branch point, leaving
the first sheet and moving to this one. 

Indeed, the Riemann hypothesis can be taken as being equivalent to
the statement that the Riemann zeroes always pass through the branch
point in moving from one sheet to another. Taking $s=\sigma+i\tau$,
if $\sigma$ is held constant at $\sigma<1/2$ while $\tau$ is increased,
the polylogarithm zeroes never pass through the branch cut, and stay
on the main sheet; whereas for $\sigma>1/2$, they do pass through
the branch cut, but not at the branch point. 
\end{figure}

\begin{figure}

\caption{\label{cap:N=3D-1 Lower-sheet}Lower sheet of the polylogarithm}

\includegraphics[%
  width=1\textwidth]{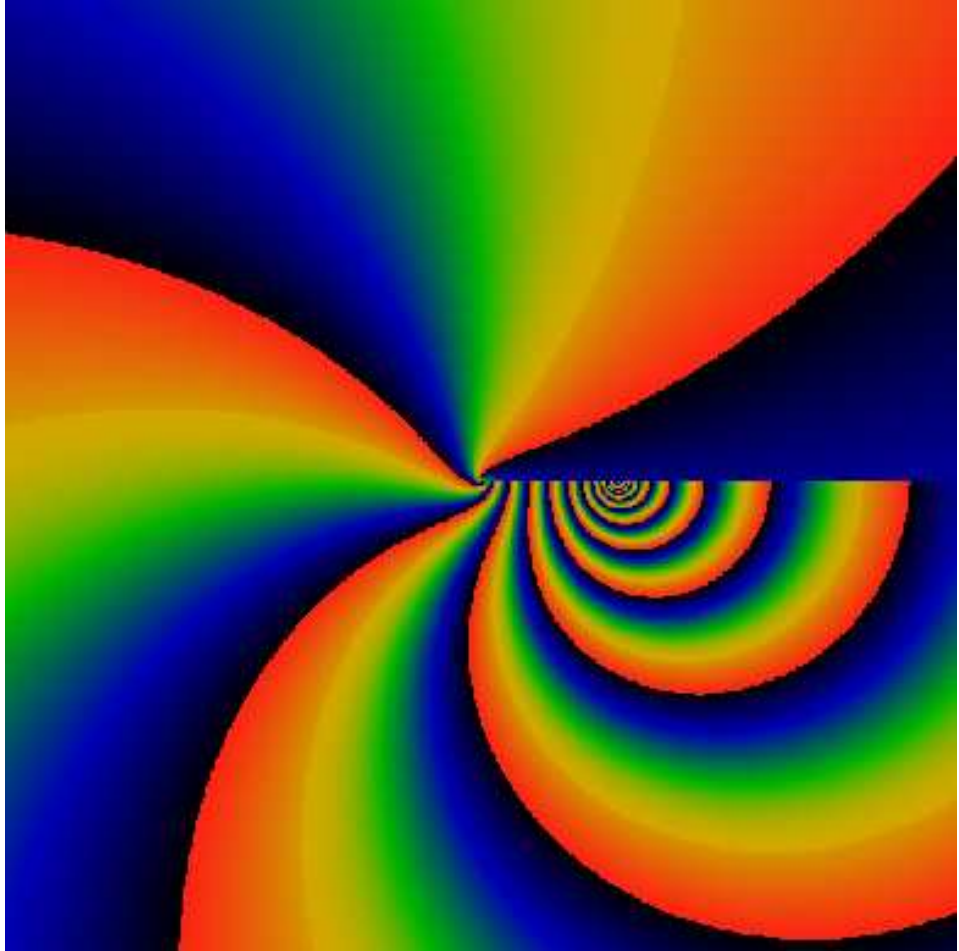}

This figure shows the $g_{1}^{-1}\mbox{Li}_{s}(z)=\mbox{Li}_{s}(z)-\Delta_{-1}$
sheet of the polylogarithm on the complex $z$-plane, obtained by
going around the branch-point at $z=1$ in the left-handed (clockwise)
direction. The concentric curves are centered on $z=1$, and whose
form is essentially that of $(1-q)^{1-s}$. The vertex of the black
triangle approaches $z=0$. By visual inspection, it is clear how
to glue this sheet to the principal sheet shown in figure \ref{cap:Polylog-whole-plane}.
After this gluing, a cut remains between the points $z=0$ and $z=1$.
This cut may be glued to the sheet that results by winding around
the branch at $z=0$ in a clockwise manner. The result is shown in
the next figure. 
\end{figure}

\begin{figure}

\caption{More polylogarithm sheets}

\includegraphics[%
  width=0.33\textwidth]{polylog-m1-sheet.ps}\includegraphics[%
  width=0.33\textwidth]{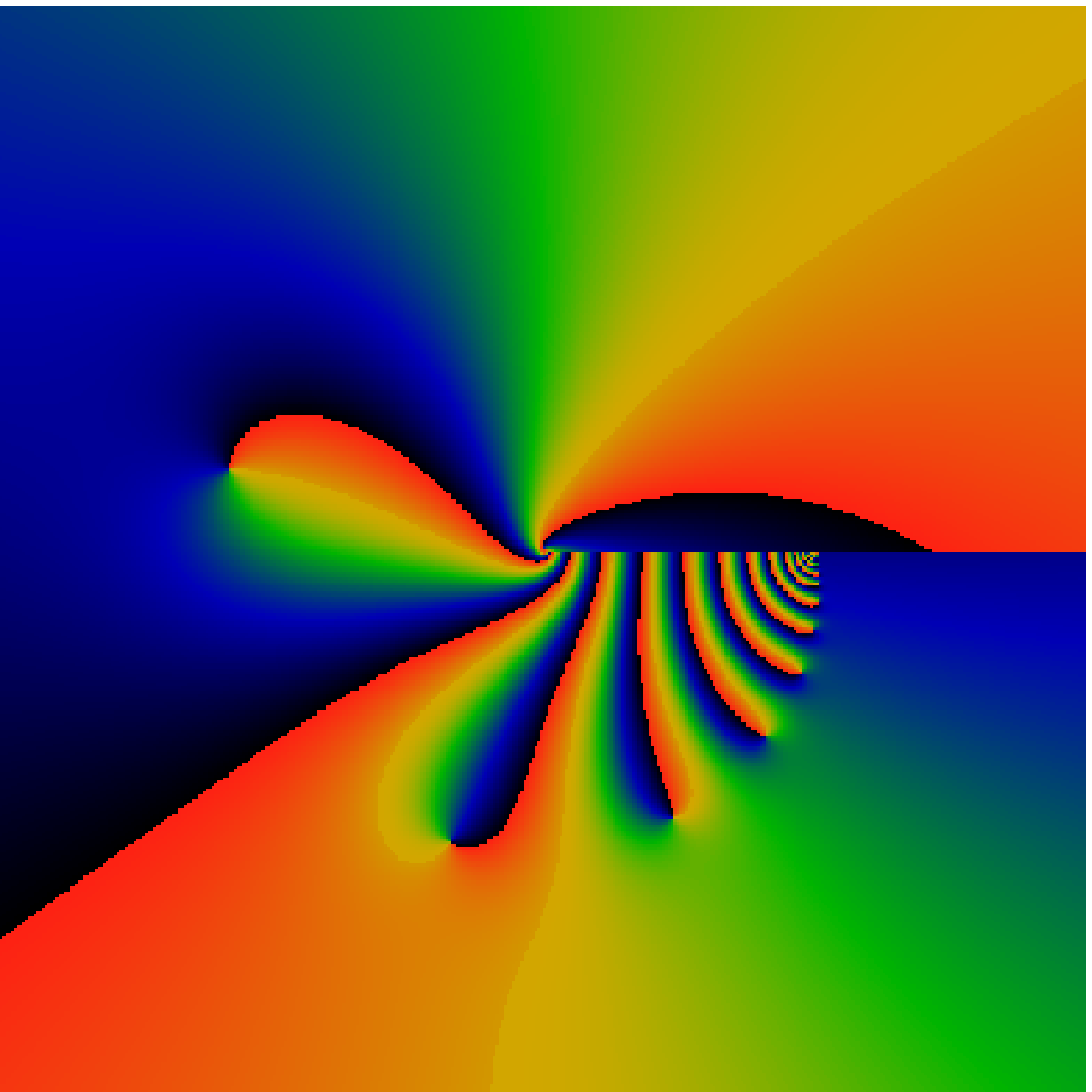}\includegraphics[%
  width=0.33\textwidth]{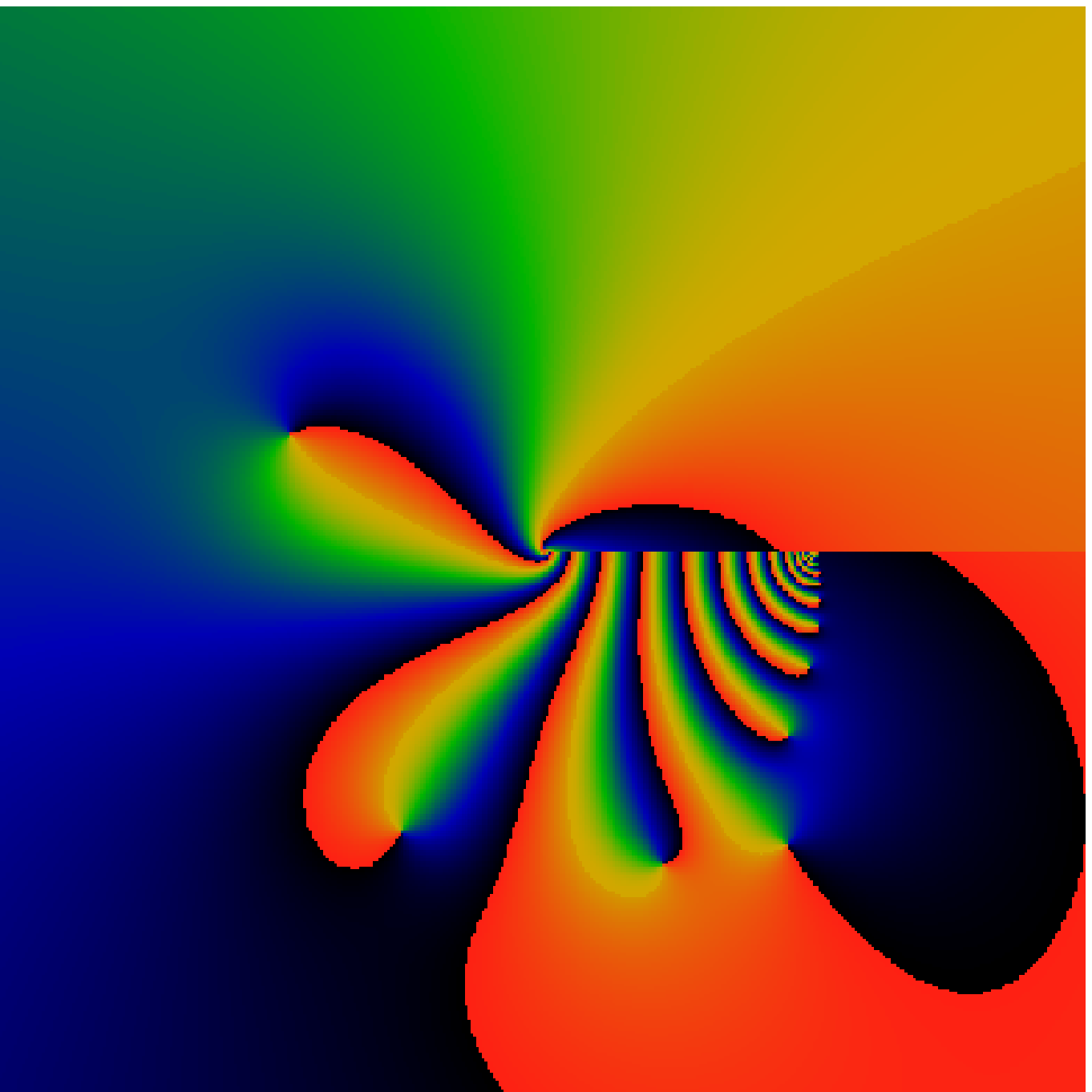}

\includegraphics[%
  width=0.33\textwidth]{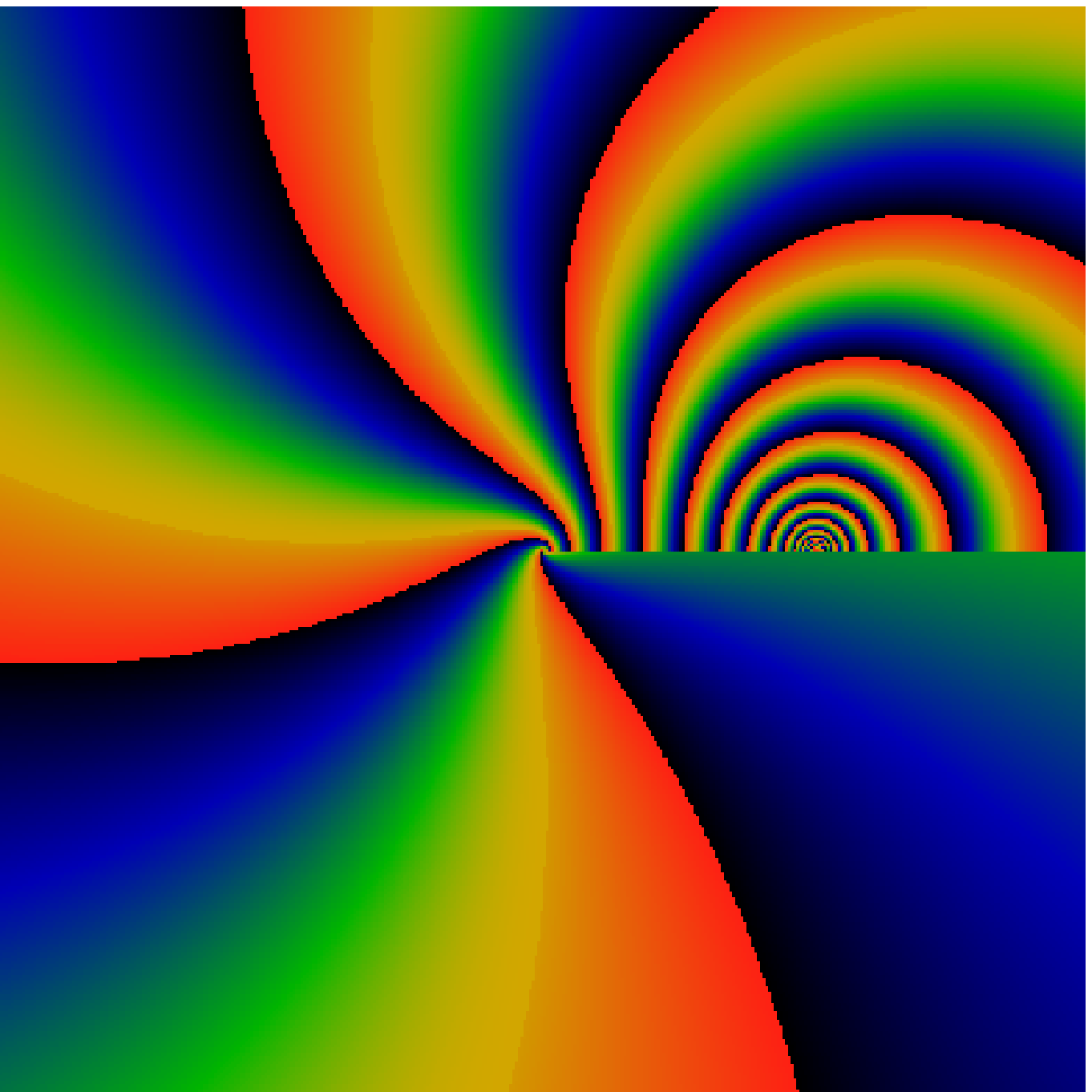}\includegraphics[%
  width=0.33\textwidth]{polylog,1,1-sheet.ps}\includegraphics[%
  width=0.33\textwidth]{polylog,1,1-sheet.ps}

Top row: shows the result of winding around $z=1$, one, two and three
times in a left-handed fashion. The cuts to the right of $z=$1 join
together smoothly from image to image. Using the monodromy group notation,
these are the $g_{1}^{-1}$, the $g_{1}^{-1}g_{0}^{-1}g_{1}^{-1}$,
and the $g_{1}^{-1}g_{0}^{-1}g_{1}^{-1}g_{0}^{-1}g_{1}^{-1}$ sheets.

Second row: Result of winding left-handed around $z=1$, followed
by winding around $z=0$ once. The cut $0<z<1$ joins smoothly to
the images in the row above. The three images appear to be visually
indistinguishable; but in fact they are not the same. Using the monodromy
group notation, these are the $g_{0}g_{1}^{-1}$, the $g_{0}g_{1}^{-1}g_{0}^{-1}g_{1}^{-1}$,
and the $g_{0}g_{1}^{-1}g_{0}^{-1}g_{1}^{-1}g_{0}^{-1}g_{1}^{-1}$
sheets. 
\end{figure}

\end{document}